\documentclass[12pt]{article}
\usepackage{amsmath}
\usepackage{amsthm}
\usepackage{graphicx}
\usepackage{psfrag}
\usepackage{amssymb}
\usepackage{amscd}
\newtheorem{theorem}{Theorem}[section]
\newtheorem{lemma}[theorem]{Lemma}
\newtheorem{corollary}[theorem]{Corollary}

\theoremstyle{definition}
\newtheorem{definition}[theorem]{Definition}

\theoremstyle{remark}
\newtheorem*{remark}{Remark}
\newtheorem*{notational remark}{Notational remark}
\newtheorem*{acknowledgements}{Acknowledgements} 
\newcommand{\Gr}{{\mathrm{Gr}}}

\newcommand{\id}{{\mathrm{id}}}
\newcommand{\Stab}{{\mathrm{Stab}}}
\newcommand{\para}{\mathrm{para}}  
\newcommand{\Area}{\mathrm{Area}}  
\newcommand{\Mod}{\mathrm{Mod}}    
\newcommand{\hol}{\mathit{hol}}
\newcommand{\QF}{\mathcal{QF}}  
\newcommand{\Q}{\mathcal{Q}}  
\newcommand{\psl}{{\mathrm{PSL}}_2(\mathbb C)}
\newcommand{\mln}{{\mathcal{ML}}_{\mathbb N}} 
\begin{document}
\title{On continuous extension of grafting maps} 
\author{Kentaro Ito}
\date{November 6, 2004}
\maketitle 
\begin{abstract}
The definition of the grafting operation 
for quasifuchsian groups   
is extended by Bromberg \cite{Br}  
to all $b$-groups. 
Although the grafting maps are not necessarily
continuous at boundary groups, 
in this paper, 
we show that the grafting maps take 
every ``standard'' convergent sequence 
to a convergent sequence. 
As a consequence of this result, 
we extend Goldman's grafting theorem 
for quasifuchsian groups to all boundary $b$-groups. 
\end{abstract}
\section{Introduction}
In this paper 
we are interested in the behavior of the holonomy map 
from projective structures to representations. 
Especially, we study the continuity and the non-continuity of the 
local inverse of the holonomy map at the boundary 
of the space of discrete faithful representations. 

Let $S$ be an oriented closed surface of genus $g>1$. 
A projective structure on $S$ is a $(G,X)$-structure 
where $X$ is a Riemann sphere $\widehat{\mathbb C}$ and 
$G=\psl$ is the group of projective automorphism of $\widehat{\mathbb C}$. 
Let $P(S)$ denote the space of projective structures 
on $S$ and  $R(S)$ the space of conjugacy classes 
of representations $\rho:\pi_1(S) \to \psl$. 
The holonomy map $\hol:P(S) \to R(S)$ 
takes a projective structure to its 
holonomy representation. 

We are interested in the quasifuchsian space 
$\QF(S) \subset R(S)$ of faithful representations with 
quasifuchsian images.  
Although $\QF=\QF(S)$ is a connected, 
contractible, open submanifold of $R(S)$,  
the topology of its closure $\overline{\QF}$ in $R(S)$ is very complicated; 
for instance McMullen \cite{Mc} showed that $\QF$ self-bumps, 
i.e., there exists $\rho \in \partial \QF$ 
such that $U \cap \QF$ is disconnected 
for any sufficiently small neighborhood $U$ of $\rho$, 
by using projective structures and ideas of Anderson and Canary \cite{AnCa}. 
In fact, since the holonomy map $\hol:P(S) \to R(S)$ 
is a local homeomorphism, 
some kind of complexity of the topology of $\QF \subset R(S)$ 
can be observed by studying 
the inverse image $Q(S)=\hol^{-1}(\QF)$ of $\QF$ in $P(S)$. 
By Goldman's grafting theorem \cite{Go}, 
the set of components of $Q(S)$ is in one-to-one correspondence 
with the set $\mln$ of integral points of measured laminations. 
For every $\lambda \in \mln$,  
the associated component $\Q_\lambda$ of $Q(S)$ 
is taken by the holonomy map biholomorphically onto $\QF$. 
We let $\Psi_{\lambda}:\QF \to \Q_\lambda$ denote 
the univalent local branch of $\hol^{-1}$ and 
call it the {\it grafting map} for $\lambda \in \mln$. 
The {\it standard component} $\Q_0$ 
is the component of $Q(S)$ of projective structures 
with injective developing map. Then 
McMullen \cite{Mc} actually showed that there exists a sequence 
in $Q(S)-\Q_0$ converging to $\partial \Q_0$, 
which implies that $\QF$ self-bumps. 
Moreover, we showed in \cite{It1,It3} that 
for any $n \in {\mathbb N}$, there exists 
$n$-components 
of $Q(S)$ which bump simultaneously 
(i.e., having intersecting closures), 
that any component of $Q(S)$ except for $\Q_0$ self-bumps, 
and that any two components of $Q(S)$ bump. 
All these phenomena are obtained by using 
exotic convergent sequence in $\QF$ (see definition below) 
constructed by Anderson and Canary \cite{AnCa}. 

Let $B:T(S) \times T(S)  \to \QF$ be the parameterization of $\QF$ 
of Bers' simultaneous uniformization and  
let $\rho_n=B(X_n,Y_n) \in \QF$ be a sequence converging 
to $\rho \in \partial \QF$. 
Then we say that the convergence $\rho_n \to \rho_\infty$ 
is {\it standard} 
if there exists a compact subset $K$ of $T(S)$ which contains 
all $X_n$ or all $Y_n$; otherwise it is {\it exotic}. 
We let $\partial^\pm \QF \subset \partial \QF$ denote the set 
of limits of standard convergence sequences, 
or the set of all boundary $b$-groups (see \S 2.1). 
Then by using Bromberg's observation in \cite{Br}, 
the grafting map $\Psi_\lambda$ is naturally 
extended to $\Psi_\lambda:\QF \sqcup \partial^\pm \QF \to \widehat{P}(S)$ 
for every $\lambda \in \mln$. 
Here $\widehat{P}(S)=P(S) \cup \{\infty\}$ 
denotes the one-point compactification of $P(S)$. 
As described above, Anderson and Canary \cite{AnCa} showed that 
there exist sequences $ \rho_n \in \QF$ converging exotically 
to $\rho_\infty \in \partial^\pm \QF$, 
which course the non-continuity of the map 
$\Psi_\lambda:\QF \to P(S)$ 
at $\rho_\infty \in \partial^\pm \QF$ 
for all $\lambda \in \mln$. 

On the contrary, 
we study the continuity of the grafting maps in this paper. 
The following theorem 
(Theorem 5.1) is the main result of the paper, 
which states that each grafting map $\Psi_\lambda$ 
behaves as a continuous map for every standard convergent sequence.
\begin{theorem}
Let 
$\rho_n \in \QF$ be a sequence converging standardly to 
$\rho_\infty \in \partial^\pm \QF$. 
Then the sequence $\Psi_{\lambda}(\rho_n)$ converges to 
$\Psi_{\lambda}(\rho_\infty)$ in $\widehat{P}(S)$ 
for every $\lambda \in \mln$. 
\end{theorem}
Let $\overline{\Q_0}$ denote the closure of $\Q_0$ in $P(S)$, 
{\it not} in $\widehat{P}(S)$.   
As we see in \S 3.3, 
a sequence $\Sigma_n \in \Q_0$ converges to $\Sigma_\infty \in \partial \Q_0$ 
if and only if 
the sequence $\hol(\Sigma_n) \in \QF$ 
converges standardly to 
$\hol(\Sigma_\infty) \in \partial^+ \QF$.  
Hence we have the next corollary (Theorem 6.10) of Theorem 1.1. 
\begin{corollary}
The map $\Gr_{\lambda}:=\Psi_\lambda \circ \Psi_0^{-1}:\Q_0 \to \Q_{\lambda}$ 
is extended continuously 
to $\Gr_\lambda:\overline{\Q_0} \to \widehat{P}(S)$ 
for each $\lambda \in \mln$, where  
$\overline{\Q_0}$ is the closure of $\Q_0$ in $P(S)$. 
\end{corollary}
An important remark is that, 
since we do not know whether $\Q_0$ self-bumps or not, 
the statement of Theorem 1.1 (or Corollary 1.2) 
is not trivial and worth considering. 

The following theorem plays an important roll 
in the proof of Theorem 1.1, 
which have some interest in its own light (see \S 4). 
A Bers slice ${\mathcal B}$ is a slice of 
$\QF$ of the type 
$\{B(X,Y)\}_{Y \in T(S)}$ or $\{B(X,Y)\}_{X \in T(S)}$. 
In addition, we let $\pi:P(S) \to T(S)$ denote 
the projection which takes a projective structure to 
its underling conformal structure. 
\begin{theorem}
Let ${\mathcal B}$ be a Bers slice and 
$\{\lambda_n\}$ a sequence of distinct elements of $\mln$. 
Then the sequence 
$\{\pi \circ \Psi_{\lambda_n}({\mathcal B})\}$ eventually escapes  
any compact subset $K$ of $T(S)$; 
that is, 
$\pi \circ \Psi_{\lambda_n}({\mathcal B}) \cap K=\emptyset$ 
for all large enough $n$. 
\end{theorem}
As a consequence of Theorem 1.1, 
we extend Goldman's grafting theorem 
as follows (see Theorem 6.1): 
\begin{theorem}
For every boundary $b$-groups $\rho \in \partial^{\pm} \QF$, all 
projective structures with holonomy $\rho$ 
are obtained by grafting of $\rho$;    
that is, 
$\hol^{-1}(\rho)=\{\Psi_\lambda(\rho) : 
\lambda \in \mln, \, \Psi_{\lambda}(\rho) \ne \infty\}.$
\end{theorem}
As mentioned above, 
any two components of $Q(S)=\hol^{-1}(\QF)$ bump. 
On the contrary, we obtain the following theorem (Theorem 6.11),
which reveals that 
only exotic convergent sequences cause the 
bumping of distinct components of $Q(S)$.  
\begin{theorem}
Set 
$\QF_K=\{B(X,Y)\in \QF:X \in K \text{ or } Y \in K \}$ 
for a compact subset $K \subset T(S)$. 
Then the inverse image 
$\hol^{-1}(\QF_K)$ of $\QF_K$ in $P(S)$ is discrete; 
that is, every connected component $\Psi_\lambda(\QF_K)$ 
of $\hol^{-1}(\QF_K)$ has an 
open neighborhood in $P(S)$ which is disjoint from any other one. 
\end{theorem}
\begin{remark}
It is conjectured by Bers, Sullivan and Thurston that 
the closure of $\QF$ is equal to 
the space $AH(S) \subset R(S)$ of discrete, faithful representations. 
This conjecture is closely related to Thurston's 
ending lamination conjecture, which was recently announced by Minsky 
to be solved affirmatively.  
But we do not make use of these deep results in this paper. 
\end{remark}
\begin{acknowledgements}
The author would like to thank Katsuhiko Matsuzaki 
for useful discussions on the topics of \S 4. 
\end{acknowledgements}

\section{Preliminaries}
\subsection{Quasifuchsian space}
A Kleinian group $\Gamma$ is a discrete subgroup 
of $\psl$, which acts on hyperbolic $3$-space ${\mathbb H}^3$ 
as isometries 
and on the sphere at infinity $\widehat{\mathbb C}$ as conformal 
automorphisms. The union 
$\overline{{\mathbb H}^3}={\mathbb H}^3 \cup \widehat {\mathbb C}$ 
is naturally topologized as a closed 3-ball 
so that $\psl$ acts continuously on it.   
For a Kleinian group $\Gamma$, 
we let $\Omega_\Gamma$ denote the region of discontinuity 
and $\Lambda_\Gamma$ the limit set. 
We associate to a Kleinian group $\Gamma$ the following orbit spaces: 
$$
M_\Gamma={\mathbb H}^3/\Gamma, \quad 
\overline{M_\Gamma}=({\mathbb H}^3 \cup \Omega_\Gamma)/\Gamma,  \quad  
 \partial M_\Gamma=\Omega_\Gamma/\Gamma,
$$ 
where $\partial M_\Gamma$ is called the 
{\it conformal boundary} of $M_\Gamma$. 
In general if $\overline{M}$ is an oriented 
maniford  with boundary $\partial M$,  
we orient $\partial M$ by requiring that the frame 
$(f,n)$ has positive orientation 
whenever $f$ is a positively oriented frame on $\partial M$ 
and $n$ is an inward-pointing vector. 

Let $S$ be an oriented closed surface of genus $g>1$. 
Let $R(S)$ be the space of conjugacy classes 
$[\rho]$ of representations 
$\rho:\pi_1(S) \to \psl$ whose images $\rho(\pi_1(S))$ are non-abelian. 
The space $R(S)$ is equipped 
with the algebraic topology, 
the topology of convergence on generators up to conjugation. 
(By abuse of notation, we also denote $[\rho]$ by $\rho$ if there is no confusion.)  
It is known that $R(S)$ is a $(6g-6)$-dimensional complex manifold 
(see Theorem 4.21 in \cite{MaTa}). 

Let $AH(S)$ be the subset of $R(S)$ 
of discrete, faithful representations, which is known to be
a closed subset of $R(S)$ by J{\o}rgensen \cite{Jo}. 
Let $\rho_n \to \rho_\infty$ be an algebraic convergence 
sequence in $AH(S)$. 
Then it is known by J{\o}rgensen and Marden \cite{JoMa} that,  
by passing to a subsequence if necessary, 
the sequence 
$\Gamma_n=\rho_n(\pi_1(S))$ 
of Kleinian groups converges geometrically to 
some Kleinian group $\widehat{\Gamma}$, which 
contains the algebraic limit 
$\Gamma_\infty=\rho_\infty(\pi_1(S))$. 
Here the convergence $\Gamma_n \to \widehat{\Gamma}$ 
is {\it geometric} if and only if  
for any $\hat{\gamma} \in \widehat{\Gamma}$, there exist 
$\gamma_n \in \Gamma_n$ such that $\gamma_n \to \hat{\gamma}$ 
and if for every convergent sequence 
$\gamma_{n_j} \in \Gamma_{n_j} \, (n_j \to \infty)$,  
the limit is contained in $\widehat{\Gamma}$. 
A convergence $\rho_n \to \rho_\infty$ in $AH(S)$ 
is said to be {\it strong} 
if $\Gamma_n=\rho_n(\pi_1(S))$ converges geometrically 
to the algebraic limit $\Gamma_\infty=\rho_\infty(\pi_1(S))$.  

For every $\rho \in AH(S)$ with image $\Gamma=\rho(\pi_1(S))$, 
Bonahon's theorem \cite{Bo} guarantees that 
there exists an orientation preserving homeomorphism 
$\psi: S \times (-1,1) \to M_\rho=M_\Gamma$, where 
$(-1,1)$ is an open interval.  
The conformal boundary $\partial M_\rho$ (possibly empty) of $M_\rho$ 
decomposes into two parts 
$\partial^+ M_\rho \sqcup \partial^- M_\rho$,  
where $\partial^+ M_\rho$ (resp. $\partial^- M_\rho$) 
is the limit of $\psi(S \times \{t\})$ in $\overline{M_\rho}$ 
as $t \to -1$ (resp. $t \to 1)$. 
Associated to this decomposition, 
$\Omega_\Gamma$ decomposes into 
$\Omega^+_\Gamma \sqcup \Omega^-_\Gamma$ 
for which 
$\Omega^+_\Gamma/\Gamma=\partial^+ M_\rho$ 
and $\Omega^-_\Gamma/\Gamma=\partial^- M_\rho$. 

Let $\Gamma=\rho(\pi_1(S))$ be the image of some $\rho \in AH(S)$.  
Then $\Gamma$ is called a {\it $b$-group} 
if exactly one of $\Omega^+_\Gamma$ or $\Omega^-_\Gamma$ 
is simply connected and $\Gamma$-invariant.  
For a $b$-group $\rho$, 
we denote by $\para(\rho)$ its {\it parabolic locus},  
a collection of homotopy classes of disjoint 
simple closed curves $c$ on $S$ 
such that $\rho(c) \in \Gamma$ is parabolic. 
Similarly $\Gamma$ is called a {\it quasifuchsian group} if 
both $\Omega^+_\Gamma$ and $\Omega^-_\Gamma$ 
are simply connected and $\Gamma$-invariant.  
Quasifuchsian space $\QF=\QF(S)$ is the subset of $R(S)$ 
of faithful representations 
with quasifuchsian images. 
It is known by Marden \cite{Mar} and Sullivan \cite{Su}
 that $\QF$ equals the interior of $AH(S)$. 
Hence $\QF$ is a $(6g-6)$-dimensional complex manifold in $R(S)$. 
On the other hand, it is trivial that $\overline{\QF} \subset AH(S)$ 
and is conjectured that $\overline{\QF} = AH(S)$,  
which is so called the Bers-Sullivan-Thurston conjecture. 
We let $\partial \QF$ denote the relative 
boundary of $\QF$ in $R(S)$, whose element is called 
a {\it boundary group}.  

Now let $\rho \in \QF$. 
Then there exists an orientation preserving homeomorphism 
$\psi: S \times [-1,1] \to \overline{M_\rho}$ 
which induces the representation 
$\psi_*=\rho:\pi_1(S) \to \psl$. 
Moreover, it provides orientation preserving homeomorphisms 
$\psi|_{S \times \{-1\}}:S \to \partial^+ M_\rho$ 
and $\psi|_{S \times \{1\}}:\bar{S} \to \partial^- M_\rho$, 
where $\bar{S}$ denotes $S$  with its orientation reversed. 
Hence $\rho \in \QF$ determines 
a pair of marked Reimann surfaces 
$(\partial^+ M_\rho, \partial^- M_\rho) \in T(S) \times T({\bar S})$ 
in the product of the Teichm\"{u}ller spaces.  
On the other hand, Bers \cite{Be} showed  
that each pair 
$(X,\bar{Y}) \in T(S) \times T(\bar{S})$ 
has the unique simultaneous uniformization 
$\rho=B(X,\bar{Y}) \in \QF$. 
Therefore the map 
$$
B: T(S) \times T(\bar{S}) \to \QF    
$$ 
gives us a global parameterization of $\QF$. 
We define {\it vertical} and {\it horizontal Bers slices} in $\QF$ by 
$B_X=\{B(X,\bar{Y}) : \bar{Y} \in T(\bar{S})\}$ and 
$B_{\bar{Y}}=\{B(X,\bar{Y}) : X\in T(S)\}$, respectively.  
It is known by Bers that both $B_X$ and $B_{\bar{Y}}$ are 
precompact in $R(S)$, whose frontiers are 
denoted by $\partial B_X$ and $\partial B_{\bar{Y}}$. 
A representation $\rho \in AH(S)$ is called a 
{\it Bers boundary group} if it is contained in 
the frontier of some Bers slice.  
  
\subsection{Sequences in quasifuchsian space}
We introduce the notion of ``standard'' and ``exotic'' 
convergence for a sequence $\rho_n \in \QF$ tending to a limit  
$\rho_\infty \in \partial \QF$.  
One of our main purpose of this paper is to show that every grafting map 
behaves as a continuous map for every 
standard convergent sequence (see Theorem 5.1). 
For a given subset $K$ of $T(S)$, we set 
$\bar{K}=\{\bar{X} \in T(\bar{S}) : X \in K\} \subset T(\bar{S})$, 
where $\bar{X} \in T(\bar{S})$ denote 
the complex conjugation of $X \in T(S)$. 
\begin{definition}[Standard and exotic convergence] 
Suppose that a sequence $\rho_n=B(X_n, \bar{Y}_n)$ in $\QF$ 
converges to $\rho_{\infty} \in \partial \QF$. 
Then the sequence $\rho_n$ is said to converge {\it standardly} 
to $\rho_\infty$ if there exists a compact subset $K \subset T(S)$ 
such that 
(i) $X_n  \in K$ for all $n$, 
or
(ii) 
$\bar{Y}_n \in \bar{K}$ for all $n$.   
Otherwise, we say that $\rho_n$ converges 
{\it exotically} to $\rho_{\infty}$. 
\end{definition}
We let $\partial^+ \QF$ and $\partial^- \QF$ 
denote the subsets of $\partial \QF$ 
consisting of limits of standard convergent sequences 
of type (i) and (ii), respectively,  
and set $\partial^\pm \QF=\partial^+ \QF \sqcup \partial^- \QF$. 
It is easily seen that every Bers boundary group 
is an element of $\partial^\pm \QF$ and that every element 
in $\partial^\pm \QF$ is a {\it boundary $b$-group}, 
a $b$-group in $\partial \QF$. 
On the contrary, it follows from the arguments in the paper of 
Brock, Bromberg, Evans and Souto \cite{BBES} 
that every boundary $b$-group is a Bers boundary group 
(see Theorem 2.3 below). 
Hence we see that the set 
$\partial^\pm \QF$ equals the set 
of all boundary $b$-groups and that the following hold: 
\begin{eqnarray*}
\partial^+\QF=\bigsqcup_{X \in T(S)}\partial B_X, \quad   
\partial^-\QF=\bigsqcup_{{\bar{Y}} \in T(\bar{S})}\partial B_{\bar{Y}}. 
\end{eqnarray*}
Moreover, Bromberg \cite{Br} obtained the following result 
by developing the grafting operation for elements 
in $\partial^- \QF$ (see \S 3) and 
by combining with Minsky's result \cite{Mi}. 
\begin{theorem}[Bromberg]
Every $b$-group $\rho \in AH(S)$ without parabolics is a Bers boundary group. 
\end{theorem}
We now outline the proof of Theorem 2.3 
by following the arguments in \cite{BBES}. 
In the argument, we also obtain Corollary 2.4 below,   
which is required in \S 5. 
\begin{theorem}[Brock-Bromberg-Evans-Souto]
Every boundary $b$-group is a Bers boundary group. 
\end{theorem}
\begin{corollary}
Let $\rho \in \partial^\pm \QF$ with parabolic locus $\para(\rho)$. 
If $\rho \in \partial^+ \QF$ (resp. $\rho \in \partial^- \QF$),   
  there exists a sequence 
  $\rho_n=B(X_n,{\bar Y}_n) \in \QF$ 
  which converges standardly to $\rho$ 
  and which satisfies 
  $l_{\bar{Y}_n}(\para(\rho)) \to 0$ 
  (resp. $l_{X_n}(\para(\rho)) \to 0$) as $n \to \infty$.  
  Here $l_{\bar{Y}_n}(\para(\rho))$ denotes 
  the total sum of hyperbolic lengths of components 
  of $\para(\rho)$ on $\bar{Y}_n$.  
\end{corollary}
\begin{proof}[Outline of the proof of Theorem 2.3]
Let $\rho \in \partial \QF$ with $\para(\rho)$, which is possibly empty. 
We may assume that the positive part $\partial^+_c M_{\rho}$  
of the conformal boundary of $M_{\rho}$ 
is homeomorphic to $S$, and hence determine a point $X \in T(S)$. 
By Theorem 3.1 in \cite{BBES}, 
the drilling theorem plays an important role in its proof, 
there exists a strong convergent sequence 
$\rho_n \to \rho$ in $\overline{\QF}$ for which 
each $\rho_n$ is a geometrically finite representation 
whose parabolic locus 
$\para(\rho_n)$ equals $\para(\rho)$. 
Since the convergence $\rho_n \to \rho$ is strong, 
$X_n =\partial^+_c M_{\rho_n}$ converges to $X$ in $T(S)$. 
Moreover, for each $n$, 
there exists a strong convergent sequence 
$\rho_{n,k} \to \rho_n \, (k \to \infty)$ 
of quasifuchsian representations 
all of whose positive parts of the conformal boundaries 
equal $X_n \in T(S)$ (see Theorem 3.4 in \cite{BBES}). 
Therefore $\rho_{n,k}=B(X_n,\bar{Y}_{n,k})$ is a sequence 
in the Bers slice $B_{X_n}$ for each $n$. 
Since the convergence $\rho_{n,k} \to \rho_n$ is strong 
and since $\rho_n$ is geometrically finite, 
the hyperbolic lengths  
$l_{\bar{Y}_{n,k}}(\para(\rho))$ 
of $\para(\rho)$ on $\bar{Y}_{n,k}$ tend to $0$ as $k \to \infty$. 
By a diagonal argument, 
we can choose a sequence 
$\rho_n'=B(X_n,\bar{Y}_n)$ from $\{\rho_{n,k}\}_{n,k \in {\mathbb N}}$
so that $\rho_n' \to \rho$ and that 
$l_{\bar{Y}_n}(\para(\rho)) \to 0$ as $n \to \infty$, 
which satisfies the desired property in Corollary 2.4. 
Now let us consider new sequence 
$\rho_n''=B(X,\bar{Y}_n)$ in $B_X$. 
Then this sequence $\rho_n'' \in  B_{X}$ 
also converges to $\rho$ 
because maximal dilatations of quasiconformal automorphism 
of $\widehat{\mathbb C}$ conjugating $\rho_n'$ to $\rho_n''$ 
tend to $1$ as $n \to \infty$. 
This implies that $\rho \in \partial B_{X}$.  
\end{proof} 
We remark that the set 
$\partial \QF-\partial^\pm \QF$ is not empty;  
for instance, it contains limits of sequences 
which appear in Thurston's double limit theorem. 
On the other hand, Anderson and Canary 
\cite{AnCa} showed that there exists a sequence in $\QF$ 
which converges exotically to some point in $\partial^\pm \QF$. 
All the known such sequences are basically obtained by their technique 
and here is a typical example: 
let $c$ be a simple closed curve on $S$ 
and let $\tau=\tau_c$ be the Dehn twist along $c$. 
Then for a fixed pair 
$(X,\bar{Y}) \in T(S) \times T(\bar{S})$, the sequence 
\begin{eqnarray*}
\rho_n=B(\tau^n \, X, \tau^{2n} \, \bar{Y})  \quad (n \in \mathbb{Z})
\end{eqnarray*} 
in $\QF$ converges exotically to some $\rho_\infty \in \partial^+ \QF$ 
as $|n| \to \infty$. 
We will make use of this sequence in \S 3.3 to explain how 
the grafting maps fail to be extended to a continuous map.

\subsection{Space of projective structures}
We only give a brief summary of projective structures 
and refer to \cite{It1} and elsewhere for more details. 

A projective structure on $S$ is a $(G,X)$-structure 
where $X$ is a Riemann sphere $\widehat{\mathbb C}$ and 
$G=\psl$ is the group of projective automorphism of $\widehat{\mathbb C}$. 
Let $P(S)$ be the space of marked projective structures on $S$. 
A projective structure $\Sigma \in P(S)$ determine the underlying 
conformal structure $\pi(\Sigma) \in T(S)$. 
It is known that 
$P(S)$ is the holomorphic affine bundle over $T(S)$ 
with the projection $\pi:P(S) \to T(S)$, and that 
$P(S)$ is a $(6g-6)$-dimensional complex manifold. 

A projective structure $\Sigma$ on $S$ 
can be lifted to that  
$\widetilde{\Sigma}$ on $\widetilde{S}$, 
where $\widetilde{S} \to S$ is the universal cover on which 
$\pi_1(S)$ acts as a covering group.  
Since $\widetilde{\Sigma}$ is simply connected, 
we obtain a developing map $f_{\Sigma}:\widetilde{S} \to \widehat{\mathbb C}$ 
by continuing charts of $\widetilde{\Sigma}$ analytically, which induce 
a holonomy representation 
$\rho_{\Sigma}: \pi_1(S) \cong \pi_1(\Sigma) \to \psl$ 
satisfying $f_\Sigma \circ \gamma=\rho_\Sigma(\gamma) \circ f_\Sigma$ for 
every $\gamma \in \pi_1(S)$. 
We remark that the pair $(f_{\Sigma},\rho_{\Sigma})$ is determined 
uniquely up to $\psl$. 
We now define the {\it holonomy map} 
$$
\hol:P(S) \to R(S)
$$ 
by $\Sigma \mapsto [\rho_{\Sigma}]$. 
Hejhal \cite{He} showed that the map $\hol$ is a local homeomorphism and 
Earle \cite{Ea} and Hubbard \cite{Hu} independently
showed that the map is holomorphic. 
\begin{theorem}[Hejhal, Earle and Hubbard]
The holonomy map $\hol:P(S) \to R(S)$ is a holomorphic local homeomorphism. 
\end{theorem}
We denote by $Q(S)=\hol^{-1}(\QF)$ 
the set of projective structures with quasifuchsian holonomy. 
An element of $Q(S)$ is said to be {\it standard} if 
its developing map is injective; otherwise it is {\it exotic}. 
We denote by  $\Q_0 \subset Q(S)$ 
the subset of standard projective structures. 
Let $\rho=B(X, \bar{Y}) \in \QF$  
with image $\Gamma=\rho(\pi_1(S))$.  
Then the quotient surface 
$\Sigma=\Omega_{\Gamma}^+/\Gamma$ can be regarded as 
a standard projective structure on $S$ 
with bijective developing map 
$f_{\Sigma}:\widetilde{S} \to \Omega_{\Gamma}^+$,  
with holonomy representation $\rho_{\Sigma}=\rho$, and with  
underlying conformal structure $X \in T(S)$. 
Let 
\begin{eqnarray*}
\Psi_0:\QF \to \Q_0
\end{eqnarray*} 
be the map defined by the correspondence 
$\rho \mapsto \Omega_{\Gamma}^+/\Gamma$ as described above. 
Then the map $\Psi_0$ turns out to be 
a univalent local branch of $\hol^{-1}$ 
onto the connected component $\Q_0$ of $Q(S)$, 
which is called the {\it standard component}. 
It is known by Bers that 
every Bers slice $B_X \subset \QF$ is embedded by the map $\Psi_0$ 
into a bounded domain $\Psi_0(B_X)$ of a fiber 
$\pi^{-1}(X) \subset P(S)$. 
Observe that $\hol|_{\Q_0}:\overline{\Q_0} \to \QF \sqcup \partial^+ \QF$ 
is bijective, where $\overline{\Q_0}$ is the closure of $\Q_0$ 
in $P(S)$.

\section{Grafting}
\subsection{Grafting maps on quasifuchsian space}
We let $\mln=\mln(S)$ denote the set of 
integral points of measured laminations on $S$. 
In other words, each element of $\lambda \in \mln$ is a 
isotopy class of disjoint union $\sqcup_{i=1}^l k_ic_i$ of 
homotopically distinct simple closed curves $c_i$ on $S$ 
with positive integer $k_i$ weights. 
We do not distinguish the isotopy class $\lambda$ and its representative 
if there is no confusion. 
The ``zero-lamination'' $0$ is also contained in $\mln$. 
In what follows, the parabolic locus $\para(\rho)$ of 
a $b$-group $\rho$ is also regarded as an element of $\mln$.  

For each non-zero $\lambda \in \mln$,  
we will explain haw to obtain the grafting map 
\begin{eqnarray*}
\Gr_\lambda:\Q_0 \to P(S),   
\end{eqnarray*} 
which satisfies $\hol \circ \Gr_\lambda \equiv \hol$ on $\Q_0$. 
We give here two equivalent definitions 
of grafting operation; 
the first one is as usual, 
and the second one is 
introduced by Bromberg in \cite{Br} so that it also 
makes sense for elements of $\partial^- \QF$.  
For a while, we assume that $\lambda$ 
is a simple closed curve $c$ of weight one for simplicity. 
In addition, we fix our notation as follows: 
for a given $\rho \in \QF$, 
let $\Gamma$ be the quasifuchsian image of $\rho$ 
and $\Omega_{\Gamma}=\Omega^+_{\Gamma} \sqcup \Omega^-_{\Gamma}$ 
the region of discontinuity of $\Gamma$.  
We regard quotient surfaces 
$\Sigma=\Omega_{\Gamma}^+/\Gamma=\Psi_0(\rho_0)$ and 
$\Sigma^-=\Omega_{\Gamma}^-/\Gamma$ as projective structures  
on $S$ and $\bar{S}$, respectively.  
Let $c^+ \subset \Sigma$ and $c^- \subset \Sigma^-$ be 
simple closed curves corresponding to $c \subset S$ 
and suppose that $\gamma \in \Gamma \cong \pi_1(S)$ is a 
representative of the homotopy class of $c$. 
Let $\tilde{c}^+ \subset \Omega^+_{\Gamma}$ and 
$\tilde{c}^- \subset \Omega^-_{\Gamma}$
be the $\langle \gamma \rangle$-invariant lifts of 
$c^+ \subset \Sigma$ and $c^- \subset \Sigma^-$, respectively.   
\begin{notational remark}
Let $A,\,B$ be open subsets of a topological space ${\mathcal X}$ 
such that $A \cap B=\emptyset$ and that $\overline{A} \cap \overline{B} \ne \emptyset$. 
Then $A \amalg B$ denotes the interior of $\overline{A \sqcup B}$.
\end{notational remark}
\begin{definition}[Grafting] 
In the situation as described above, 
the {\it grafting}  $\Gr_{\lambda}(\Sigma)$ 
of the standard projective structure 
$\Sigma=\Psi_0(\rho)$ along $c$ is a projective structure 
obtained by the following (equivalent) procedures: 
\begin{description}
  \item[I] 
let $A_c$ be a cylinder 
$(\widehat{\mathbb C}-\tilde{c}^+)/\langle \gamma \rangle$ 
equipped with a projective structure 
induced from that of $\widehat{\mathbb C}$. 
We obtain $\Gr_c(\Sigma)$ by 
cutting $\Sigma$ along $c$ and inserting 
$A_c$ at the cut locus without twisting; that is, 
$$\Gr_c(\Sigma)=(\Sigma-c) \amalg A_c.$$
  \item[II]   
Here we further assume that 
$c$ separates $S$ into two surfaces  
$S_1$ and $S_2$ with boundaries. 
(The non-separating case is described precisely in \cite{Br}.) 
Accordingly, $\Sigma$ and $\Sigma^-$ decompose into 
$\Sigma-c^+=\Sigma_1 \sqcup \Sigma_2$ and 
$\Sigma^- -c^-=\Sigma^-_1 \sqcup \Sigma^-_2$, respectively. 
Let $i$ denotes $1$ or $2$. 
Let $\Delta_i$ be the connected component of the inverse image of 
$\Sigma_i^- \subset \Sigma^-$ in $\Omega_{\Gamma}^-$ 
which contains $\tilde{c}^-$  
in its closure $\overline{\Delta}_i \subset \widehat{\mathbb C}$. 
Then the stabilizer subgroup 
$\Gamma_i=\Stab_{\Gamma}(\Delta_i)$ of $\Gamma \cong \pi_1(S)$ 
is identified with $\pi_1(S_i)$. 
Since $\Gamma_i$ is a purely loxodromic free group 
with non-empty region of discontinuity, 
Maskit's result \cite{Mas} implies that $\Gamma_i$ is a Schottky group. 
Note that the conformal boundary 
$\partial M_{\Gamma_i}=\Omega_{\Gamma_i}/\Gamma_i$ of $M_{\Gamma_i}$ 
with natural projective structure containing 
both projective surfaces $\Sigma_i$ and $\Sigma_i^-$. 
Then $\Gr_{c}(\Sigma)$ is obtained from 
projective surfaces 
$\Omega_{\Gamma_1}/\Gamma_1-\Sigma_1^-$ and 
$\Omega_{\Gamma_2}/\Gamma_2-\Sigma_2^-$ 
by gluing their boundaries without twisting
(see Figure 1); that is, 
$$\Gr_c(\Sigma)=(\Omega_{\Gamma_1}/\Gamma_1-\Sigma_1^-) 
\amalg (\Omega_{\Gamma_2}/\Gamma_2-\Sigma_2^-).$$
\end{description}
\end{definition}
\begin{figure}[htbp]
  \begin{center}
      \begin{psfrags}
\psfrag{S1+}{$\Sigma_1$}
\psfrag{S1-}{$\Sigma^-_1$}
\psfrag{S2+}{$\Sigma_2$}
\psfrag{S2-}{$\Sigma^-_2$}
\psfrag{G}{$\Gr_c(\Sigma)$}
\psfrag{M1}{$M_{\Gamma_1}$}
\psfrag{M2}{$M_{\Gamma_2}$}
 \includegraphics{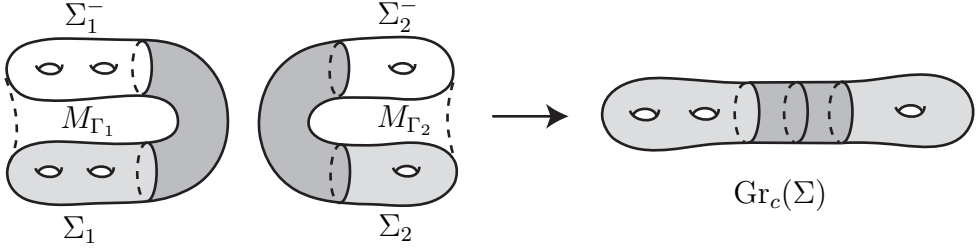}
      \end{psfrags}
     \end{center}
       \caption{The grafting $\Gr_c(\Sigma)$ of $\Sigma$ along $c$.}
\end{figure}
Observe that the Definitions I and II are equivalent.  
The grafting $\Gr_\lambda(\Sigma)$ of $\Sigma$ along general 
$\lambda=\sqcup k_ic_i \in \mln$ is similarly obtained; 
for instance, $\Gr_{kc}(\Sigma)$ is obtained by inserting 
$k$-copies of $A_c$ in Definition I. 
An important fact is that the grafting operation does not 
change the holonomy representation, that is, 
$\hol(\Gr_\lambda(\Sigma))=\hol(\Sigma)$ is always satisfied. 

Since the grafting map $\Gr_\lambda:\Q_0 \to P(S)$ satisfies 
$\hol \circ \Gr_\lambda \equiv \hol$ on $\Q_0$ and since 
$\hol|_{\Q_0}:\Q_0 \to \QF$ is a biholomorphic map with 
its inverse $\Psi_0$, 
the map $\Gr_\lambda$ takes $\Q_0$ 
biholomorphically onto the image $\Gr_{\lambda}(\Q_0)$,  
which is denoted by $\Q_\lambda$.  
Hence we obtain a univalent local branch 
$$
\Psi_\lambda:\QF \to \Q_{\lambda}
$$
of $\hol^{-1}$ which is defined by $\Psi_\lambda=\Gr_\lambda \circ \Psi_0$
and is also called the {\it grafting map} for 
$\lambda$ (see the commutative diagram below). 
\begin{figure}[htbp]
  \begin{center}
      \begin{psfrags}
\psfrag{P}{$P(S) \supset \Q_0$}
\psfrag{R}{$R(S) \supset \QF$}
\psfrag{L}{$\Q_{\lambda}$}
\psfrag{h}{$\hol$}
\psfrag{a}{$\Psi_0$}
\psfrag{b}{$\Psi_{\lambda}$}
\psfrag{g}{$\Gr_{\lambda}$}
 \includegraphics{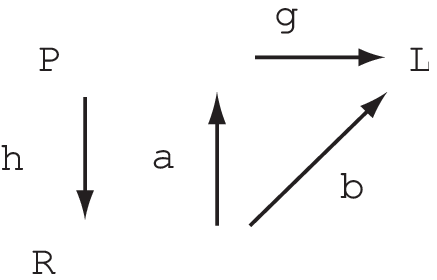}
      \end{psfrags}
     \end{center}
\end{figure}
By Goldman's grafting theorem \cite{Go} below, 
we obtain the decomposition 
$\bigsqcup_{\lambda \in \mln}{\mathcal Q}_{\lambda}$
of $Q(S)$ into its connected components. 
\begin{theorem}[Goldmann \cite{Go}]
For every $\rho \in \QF$, all projective structures 
with holonomy $\rho$ are obtained by grafting of $\rho$; 
that is, we have 
$\hol^{-1}(\rho)=\{\Psi_{\lambda}(\rho) \,:\, \lambda \in \mln\}$. 
\end{theorem}

\subsection{Extension of grafting maps}
Let $\widehat{P}(S)$ denotes the 
one-point compactification $P(S) \cup \{\infty\}$ of $P(S)$.  
We now extend the grafting map $\Psi_\lambda:\QF \to \Q_\lambda$ 
to $\Psi:\QF \sqcup \partial^\pm \QF \to \widehat{P}(S)$.  
(In fact $\Psi_\lambda(\rho)$ is defined 
for all $b$-groups $\rho \in AH(S)$ by the same manner.)   
Observe that Definition I works well even 
for $\rho \in \partial^+ \QF$ whenever $\gamma$ is loxodromic, 
because there still exists a $\langle \gamma \rangle$-invariant 
simple arc $\tilde{c}^+$ 
in non-degenerate component $\Omega^+_{\Gamma}$,  
for which $(\widehat{\mathbb C}-\tilde{c}^+)/\langle \gamma \rangle$ 
is still an annulus. 
On the other hand, Definition II works well for $\rho \in \partial^-\QF$ 
whenever every connected 
component of $\para(\rho)$ intersects $c$ essentially. 
In fact, on this assumption,  
$\Gamma_1$ and $\Gamma_2$ in Definition II are still Schottky groups,  
$\gamma$ is still loxodromic,  
and there still exists a 
$\langle \gamma \rangle$-invariant simple arc $\tilde{c}^-$ 
in non-degenerate component $\Omega^-_{\Gamma}$.  
For general $\lambda \in \mln$, 
we can also define the grafting 
$\Psi_{\lambda}(\rho) \in P(S)$ of $\rho$ along $\lambda$ 
if the pair $(\lambda,\rho)$ satisfies the following condition: 
\begin{definition}[Admissible]
The pair $(\lambda,\rho)$ of $\lambda \in \mln$ 
and $\rho \in \partial^\pm \QF$ 
is said to be {\it admissible} if 
\begin{itemize}
  \item $\rho \in \partial^+ \QF$ and 
  $\para(\rho)$ and $\lambda$ have no parallel component in common, or   
  \item $\rho \in \partial^- \QF$ and 
every component of $\para(\rho)$ 
intersects $\lambda$ essentially. 
\end{itemize}
\end{definition}
If the pair $(\lambda,\rho)$ is not admissible, 
we set $\Psi_{\lambda}(\rho)=\infty \in \widehat{P}(S)$. 
The extended grafting map 
$\Psi_\lambda:\QF \sqcup \partial^\pm \QF \to \widehat{P}(S)$ 
is also denoted by the same symbol $\Psi_\lambda$. 

\subsection{Non-continuity of grafting maps}
We collect in Table 1 below the equivalent conditions 
with standard/exotic convergence of quasifuchsian representations 
(see \cite[Appendix A]{Mc} and \cite[\S 3]{It1}). 
The situation in which we consider is as follows: 
let $\rho_n \in \QF$ be a sequence 
converging  to $\rho_{\infty} \in \partial^+ \QF$.  
We may assume that the sequence 
$\Gamma_n=\rho_n(\pi_1(S))$ 
converges geometrically to a Kleinian group $\widehat{\Gamma}$, 
which contains the algebraic limit $\Gamma_\infty=\rho_\infty(\pi_1(S))$.  
Let $\Sigma_\infty$ be a unique projective structure in $\partial \Q_0$ 
such that $\hol(\Sigma_\infty)=\rho_\infty$ and 
let $\Phi:U \to P(S)$ be a univalent local branch of 
$\hol^{-1}$ which is defined on a neighborhood $U$ of $\rho_{\infty}$ 
and takes  $\rho_{\infty}$ to $\Sigma_\infty$.  
Then $\Sigma_n=\Phi(\rho_n)$ converges to $\Sigma_\infty=\Phi(\rho_\infty)$. 
In this situation, all conditions in the same line in Table 1 are 
equivalent. 
\begin{table}[htbp]
\caption{Equivalent conditions 
with standard/exotic convergence.}
 \begin{center}
 {\tabcolsep=3mm 
 \renewcommand\arraystretch{1.5} 
  \begin{tabular}{|c||c|}
    \hline
     $\rho_n \to \rho_\infty$ : standard  &  
     $\rho_n \to \rho_\infty$ : exotic  \\
    \hline
     $\Gamma_\infty$ is a component subgroup of $\widehat{\Gamma}$  & 
     $\Gamma_\infty$ is not a component subgroup of $\widehat{\Gamma}$   \\
    \hline
     $\Omega^+_{\Gamma_\infty} \cap \Lambda(\widehat{\Gamma})=\emptyset$ &  
     $\Omega^+_{\Gamma_\infty} \cap \Lambda(\widehat{\Gamma}) \ne \emptyset$  \\
    \hline
     $\Sigma_n$ are standard \ ($n \gg 0$) & 
     $\Sigma_n$ are exotic \ ($n \gg 0$) \\
    \hline
  \end{tabular}}
 \end{center}
\end{table}

In the same setting as above, 
let $\rho_n=B(\tau^n X, \tau^{2n} \bar{Y})$ be the sequence 
introduced at the end of \S 2.1 
which converges exotically to $\rho_\infty$. 
Then $\Sigma_n$ are exotic projective structure for all 
large enough $|n|$. 
Since $\Sigma_n \to \Sigma_\infty \in \partial \Q_0$, 
the sequence $\Psi_0(\rho_n) \in \Q_0$ of standard projective structure 
can not accumulates on a point in $\partial \Q_0$ and hence diverges. 
This implies that the map 
$\Psi_0:\QF \sqcup \partial^\pm \QF \to \widehat{P}(S)$ 
is not continuous at $\rho_\infty \in \partial^+ \QF$. 
In addition, we have shown in \cite{It1} 
that $\Sigma_n=\Psi_c(\rho_n) \in \Q_c$ and that 
$\lim_{n \to +\infty}\Psi_c(\rho_n)=\lim_{n \to -\infty}\Psi_c(\rho_n)$.  
On the other hand, we have shown in \cite{It3} that 
$\lim_{n \to +\infty}\Psi_{c'}(\rho_n) 
\ne \lim_{n \to -\infty}\Psi_{c'}(\rho_n)$ 
for every simple closed curve $c'$ on $S$ which 
intersects $c$ essentially.
This implies that the grafting map 
$\Psi_{c'}:\QF \sqcup \partial^\pm \QF \to \widehat{P}(S)$ 
is not continuous at $\rho_\infty  \in \partial^+ \QF$. 
More general, we have the following 
\begin{theorem}[\cite{It3}]
For any $\lambda \in \mln$, the grafting map 
$\Psi_\lambda:\QF \sqcup \partial^\pm \QF \to \widehat{P}(S)$ 
is not continuous.  
\end{theorem}

\subsection{Pull-backs of region of discontinuities}
Let $\Sigma=(f_\Sigma,\rho_\Sigma) \in P(S)$ 
with $\Gamma=\rho_{\Sigma}(\pi_1(S))$ and let 
$X$  be a subset of $\widehat{\mathbb C}$ which is $\Gamma$-invariant. 
Then we naturally obtain a subset $f^{-1}_\Sigma(X)/\pi_1(S)$ of $\Sigma$ 
since the inverse image $f^{-1}_{\Sigma}(X) \subset \tilde{\Sigma}$ of $X$ 
is also invariant under the action of the covering group 
$\pi_1(S)$ of the covering map $\widetilde{\Sigma} \to \Sigma$. 
We call $f^{-1}_{\Sigma}(X)/\pi_1(S)$ the {\it pull-back} of $X$ in $\Sigma$. 

For $\lambda = \sqcup_i k_ic_i \in \mln$, 
a {\it realization} $\widehat{\lambda}$ of $\lambda$ is 
a disjoint union of simple closed curves on $S$ 
which realize each weighted simple closed curve $k_ic_i$ 
by $k_i$ parallel disjoint simple closed curves 
which are homotopic to $c_i$. 
Let ${\mathcal N}(\widehat{\lambda})$ denotes a regular 
neighborhood of $\widehat{\lambda}$ in $S$. 

Let $\rho \in \QF$ with $\Gamma=\rho(\pi_1(S))$, 
let $\Sigma=\Psi_\lambda(\rho)$ be the grafting of $\rho$ along $\lambda \in \mln$, 
and let $s$ denotes either $+$ or $-$.  
Then we denote by $\Omega^s_\Sigma$ the pull-back 
$f^{-1}_\Sigma(\Omega^s_\Gamma)/\pi_1(S)  \subset \Sigma$ of an 
invariant component $\Omega^s_\Gamma$ of $\Omega_\Gamma$.  
From the definition of grafting, we have the following lemma 
(see also the proof of Theorem 4.2). 
\begin{lemma}
Let $\Sigma \in \Q_\lambda$. Then 
there is a automorphism of $\Sigma$ which is homotopic to 
the identity and 
which takes $\Sigma-{\mathcal N}(\widehat{\lambda})$ onto $\Omega^+_\Sigma$ 
and ${\mathcal N}(\widehat{\lambda})$ onto $\Omega^-_\Sigma$. 
\end{lemma}
Therefore the topological type of the pull-backs 
$\Omega^+_\Sigma$ and $\Omega^-_\Sigma$ are constant  
for every elements $\Sigma \in \Q_\lambda$. 
The following lemma states that the similar result holds 
for a grafting $\Sigma=\Psi_\lambda(\rho)$ of $\rho \in \partial^\pm \QF$ 
and for the pull-back of the 
unique invariant component of $\Omega_\Gamma$ (see Figure 2). 
In fact, we use this fact in \S 5 to show that 
$\Sigma=\Psi_\lambda(\rho)$ for $\rho \in \partial^\pm \QF$ is actually a limit of 
a sequence in $\Q_\lambda$. 
Similar arguments can be found in \cite{It1} and \cite{It3}.
\begin{lemma}
Let $(\lambda,\rho)$ be an admissible pair of 
$\lambda \in \mln$ and $\rho \in \partial^\pm \QF$, and 
set $\Sigma=\Psi_\lambda(\rho)$. 
\begin{itemize}
    \item If $\rho \in \partial^+ \QF$, 
    there is a one-to-one correspondence 
    between the set of connected components 
    of $\Omega^+_\Sigma \subset \Sigma$ and that of 
    $\Sigma-{\mathcal N}(\widehat{\lambda}) \subset \Sigma$. 
    Moreover, 
      for each connected component of 
    $\omega \subset  \Omega^+_\Sigma$, 
    there is a homeomorphism from $\omega$ 
    onto the corresponding component of 
    $\Sigma-{\mathcal N}(\widehat{\lambda})$ 
    which is isotopic to the identity.                                
    \item If $\rho \in \partial^- \QF$, 
    there is a one-to-one correspondence 
    between the set of connected components 
    of $\Omega^-_\Sigma \subset \Sigma$ and that of 
    ${\mathcal N}(\widehat{\lambda}) \subset \Sigma$. 
    Moreover, 
    for each connected component of 
    $\omega \subset  \Omega^-_\Sigma$, 
    there is a homeomorphism from $\omega$ 
    onto the corresponding component of 
    ${\mathcal N}(\widehat{\lambda})$ 
    which is isotopic to the identity.   
    \end{itemize}
\end{lemma}
\begin{proof} 
We discuss in the same setting of Definition 3.1 and \S 3.2, 
except for one change; 
here $\Sigma$ denotes the grafting $\Psi_\lambda(\rho)$ but not a standard 
projective structure. 
We first suppose that $\rho \in \partial^+\QF$. 
Let us consider the quotient torus 
$T_c=\widehat{\mathbb C}/\langle \gamma \rangle$, 
in which the arc $\tilde{c}^+ \subset \Omega^+_\Gamma$ descends 
to an essential curve $c^+$ in $T_c$. 
Then the grafting annulus 
$A_c=({\widehat{\mathbb C}}-\tilde{c}^+)/\langle \gamma \rangle$ 
is identified with the subset $T_c -c^+$ of $T_c$.  
On the other hand, since the unique invariant component 
$\Omega^+_{\Gamma}$ of $\Omega_\Gamma$ is simply connected, 
$\Omega^+_{\Gamma}/\langle \gamma \rangle \subset T_c$ 
is also an annulus with core curve $c^+$.    
Hence the intersection 
$A_c \cap (\Omega^+_{\Gamma}/\langle \gamma \rangle)$ 
equals $\Omega^+_{\Gamma}/\langle \gamma \rangle-c^+$ in $T_c$,  
which is a disjoint union of two essential annuli $A_1$ and $A_2$. 
We now regard that $A_1 \sqcup A_2 \subset A_c \subset \Sigma$. 
Then the pull-back $\Omega^+_\Sigma \subset \Sigma$ 
of $\Omega^+_\Gamma$ coincides with 
the union $(\Sigma-A_c) \amalg A_1 \amalg A_2 \subset \Sigma$, 
which is obtained by gluing those pieces along copies of $c^+$. 
Thus the result follows in this case.  

Next we suppose that $\rho_\infty \in \partial^- \QF$. 
By the same argument as above, we see that the quotient 
$\Omega^-_{\Gamma}/\langle \gamma \rangle$ of the unique 
invariant component $\Omega^-_\Gamma$ of $\Omega_\Gamma$ 
is an annulus with core curve 
$c^-=\tilde{c}^-/\langle \gamma \rangle$, and that 
the disjoint union 
$\Omega^-_\Gamma-\tilde{c}^-=\Omega_1^- \sqcup \Omega_2^-$ 
descends to that of essential annuli 
$A^-_1 \sqcup A^-_2$ 
in $\Omega^-_{\Gamma}/\langle \gamma \rangle$, 
where $A^-_i=\Omega^-_i/\langle \gamma \rangle$ for $i=1,2$. 
We may assume that $\Delta_i \subset \Omega_i^-$ for $i=1,2$. 
Recall that $\Gamma_i=\Stab_\Gamma(\Delta_i)$ is a Schottky group.  
Since $\Stab_{\Gamma_1}(\Omega_2^-)=\langle \gamma \rangle$, 
the image of $\Omega_\Gamma^-=\Omega_1^- \amalg \Omega_2^-$ 
via the covering map $\Omega_{\Gamma_1} \to \Omega_{\Gamma_1}/\Gamma_1$ 
is $\Sigma^-_1 \amalg A_2^- \subset \Omega_{\Gamma_1}/\Gamma_1$, 
which is obtained from $\Sigma_1^-=\Delta_1/\Gamma_1$ 
and $A_2^-=\Omega^-_2/\langle \gamma \rangle$ by gluing along $c^-$. 
Similarly, 
the image of $\Omega_\Gamma^-=\Omega_1^- \amalg \Omega_2^-$ 
via the covering map $\Omega_{\Gamma_2} \to \Omega_{\Gamma_2}/\Gamma_2$ 
is $\Sigma^-_2 \amalg A_1^- \subset \Omega_{\Gamma_2}/\Gamma_2$. 
Since 
$\Sigma$ is obtained from 
$\Omega_{\Gamma_1}/\Gamma_1-\Sigma_1^-$ and 
$\Omega_{\Gamma_2}/\Gamma_2-\Sigma_2^-$ by gluing along $c^-$,   
we see that the pull-back 
$\Omega^-_{\Sigma} \subset \Sigma$ of $\Omega^-_\Gamma$ 
equals the annulus 
$A_1^-  \amalg A_2^-=\Omega^-_1/\langle \gamma \rangle
\amalg \Omega^-_2/\langle \gamma  \rangle
=\Omega^-_{\Gamma}/\langle \gamma  \rangle \subset \Sigma$ 
with core curve $c^-$.   
Thus we have obtained the result also in this case. 
\end{proof}
\begin{figure}[htbp]
  \begin{center}
      \begin{psfrags}
\psfrag{+}{$\Omega^+_\Sigma$}
\psfrag{-}{$\Omega^-_\Sigma$}
 \includegraphics{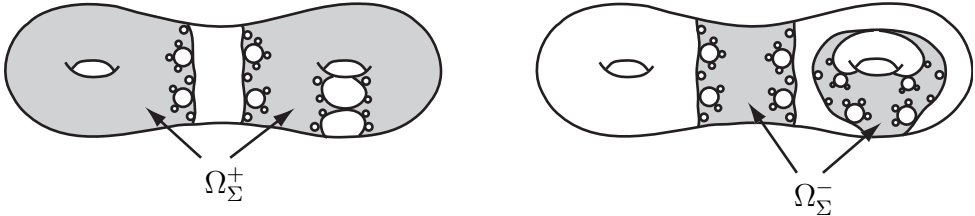}
      \end{psfrags}
     \end{center}
  \caption{Pull-backs 
  $\Omega^+_{\Sigma} \subset \Sigma$ for $\rho \in \partial^+ \QF$ and 
  $\Omega^-_{\Sigma} \subset \Sigma$ for $\rho \in \partial^- \QF$.}
\end{figure} 
\begin{corollary}
Let $\lambda,\mu\in \mln$, $\lambda \ne \mu$  
and $\rho \in \partial^\pm \QF$.  
If both pairs $(\lambda,\rho)$ and $(\mu,\rho)$ are admissible 
then $\Psi_\lambda(\rho) \ne \Psi_\mu(\rho)$.  
\end{corollary}

\section{Finiteness}
The aim of this section is to prove Theorem 4.3, 
which guarantees some finiteness:  
suppose that a sequence $\rho_n \in \QF$ converges standardly to 
$\rho \in \partial^{\pm}\QF$ and that 
$\Phi:U \to P(S)$ is a univalent local branch of $\hol^{-1}$ defined on 
a neighborhood $U$ of $\rho$. 
Then the sequence $\Phi(\rho_n)$ 
is actually contained in some finite union of components of $Q(S)$. 
This finiteness theorem plays an important role 
in the proof of our main theorem (Theorem 5.1), 
which guarantees some continuity of the grafting maps and 
which can be seen also as the uniqueness theorem. 

\subsection{Length-Modulus inequality}
Recall that the {\it modulus} $\Mod(A)$ of 
a conformal annulus $A$ is defined uniquely 
as the ratio of the height and the circumference of an 
Euclidian annulus which is conformally equivalent to $A$.  
The following inequality is a direct consequence 
of the geometric and analytic definitions of the extremal length.  
\begin{lemma}
Let $R$ be a complete hyperbolic surface of finite area 
$\Area(R)<\infty$ and let $A$ be an essential annular domain in $R$.  
Then we have 
\begin{eqnarray*}
\Mod(A)\, l_R(c)^2 \le \Area(R), 
\end{eqnarray*}
where $l_R(c)$ is the hyperbolic length of 
the homotopy class $c$ of a core curve of $A$. 
\end{lemma}
\begin{proof}
Let $E_R(c)$ denotes the extremal length of the homotopy class $c$ in $R$. 
From the analytical definition of $E_R(c)$, we have 
\begin{eqnarray*}
E_R(c):=\sup_{\rho}\frac{\left( \inf_{c'} 
\int_{c'} \rho(z) \,|dz|  \right)^2}
{\iint \rho(z)^2\,|dz|^2}
 \ge \frac{(l_R(c))^2}{\Area(R)}, 
\end{eqnarray*}
where the supremum is taken over all metrics $\rho$ consistent with 
the conformal structure of $R$ and the infimum 
is taken over all closed curves $c'$ in the 
homotopy class $c$.  
On the other hand, from the geometrical 
definition of $E_R(c)$, we have 
\begin{eqnarray*}
E_R(c):=\frac{1}{\sup_{A' \subset R} \Mod(A')} \le \frac{1}{\Mod(A)}, 
\end{eqnarray*}
where the supremum is taken over all annuli 
$A' \subset R$ whose core curve is in the homotopy class $c$. 
From the above two inequality, the desired inequality follows. 
\end{proof}

\subsection{Quasiconformal deformations}
We introduce the notion of a quasiconformal 
deformation of  a projective structure with quasifuchsian holonomy, 
which was developed by Shiga and Tanigawa in \cite{ShTa}. 
Let $\Sigma \in Q(S)$ and let $\rho_\Sigma \in \QF$ be its holonomy.  
Suppose that $\rho' \in \QF$ is a quasiconformal deformation 
of $\rho_\Sigma$ induced by a quasiconformal automorphism 
$q:\widehat{\mathbb C} \to \widehat{\mathbb C}$, whose Beltrami 
differential is denoted by $\mu$.  
Then we obtain a new projective structures 
$\Sigma'$ with holonomy $\rho'$ which is characterized as follows: 
\begin{enumerate}
  \item 
There is a quasiconformal map 
$\tilde{\varphi}:\widetilde{\Sigma} \to \widetilde{\Sigma'}$ 
whose Beltrami differential is equal to 
the pull-back $f_\Sigma^*(\mu)$ of $\mu$ via the developing map 
$f_\Sigma:\widetilde{\Sigma} \to \widehat{\mathbb C}$. 
Moreover, the map 
$\tilde{\varphi}:\widetilde{\Sigma} \to \widetilde{\Sigma'}$ 
descends to a quasiconformal map 
$\varphi:\Sigma \to \Sigma'$, 
which is consistent with their markings. 
  \item 
The developing map of $\Sigma'$ is defined by 
$f_{\Sigma'}=q \circ f_\Sigma \circ \tilde{\varphi}^{-1}:
\widetilde{\Sigma'} \to \widehat{\mathbb C}$. 
\end{enumerate}
Here we say that a map between projective structures 
is {\it quasiconformal} if it is a quasiconformal map between 
their underling conformal structures.  
We call $\Sigma'$ as the {\it quasiconformal deformation} of $\Sigma$.  
We remark that every grafting map $\Psi_\lambda:\QF \to \Q_\lambda$ 
is obtained by quasiconformal deformations of some fixed 
$\Sigma \in \Q_\lambda$ and its holonomy $\rho_\Sigma \in \QF$.   

\subsection{Length estimates}
\begin{theorem}
Fix $X \in T(S)$ arbitrarily. 
Suppose that $\lambda \in \mln$ contains a weighted simple closed curve 
$kc$ of weight $k \in {\mathbb N}$.  
Let $\rho \in B_X \cup B_{\bar{X}}$ and 
let $X' \in T(S)$ denotes the underlying conformal structure 
of $\Psi_\lambda(\rho)$. 
Then we have $l_{X'}(c)^2 \le 4(g-1)\,l_X(c)/k$, 
where $g$ denotes the genus of $S$. 
\end{theorem}
\begin{remark}
We make use of the inequality above mostly in the following form: 
$$
\frac{l_{X'}(c)}{l_X(c)} \le \sqrt{\frac{4(g-1)}{k \cdot l_X(c)}}, 
$$
which ensures that if $k \cdot l_X(c)$ is large enough then 
the Teichm\"{u}ller distance $d_{T(S)}(X,X')$ 
between $X$ and $X'$ is also large. 
\end{remark}
\begin{proof}[Proof of Theorem 4.2]
We first treat the case where $\lambda$ is a simple closed curve $c$ 
of weight one.
Let $\rho^0=B(X, \bar{X}) \in \QF$ be the Fuchsian representation with 
image $\Gamma^0=\rho^0(\pi_1(S))$.   
We normalize $\rho^0$ so that $\Omega^+_{\Gamma^0}$ 
equals the upper half plane $H=\{x+iy \in {\mathbb C} : y>0\}$, 
$\Omega^-_{\Gamma^0}$ equals the lower half plane 
$L=\{x+iy \in {\mathbb C} : y<0\}$, and that the hyperbolic element 
$\gamma=\rho^0(c) \in \Gamma^0$ fixes 
the positive imaginary axis $i{\mathbb R}_+$. 
We let $\Sigma^0$ denote the standard projective structure 
$\Psi_0(\rho^0)=\Omega^+_{\Gamma^0}/\Gamma^0$. 
Recall that the grafting 
$\Psi_c(\rho^0)=\Gr_c(\Sigma^0)$ of $\Sigma_0$, 
which is denoted by $\Sigma^0_c$, is obtained from 
$\Sigma^0$ by cutting along $c$ and inserting 
the annulus $A=({\mathbb C}- i{\mathbb R}_+)/\langle \gamma \rangle$. 
We set $H-i{\mathbb R}_+=H_L \sqcup H_R$. 
Then ${\mathbb C}- i{\mathbb R}_+=H_R \sqcup L \sqcup H_L$ 
and then $A=A_1 \amalg A_2 \amalg A_3$,   
where $A_1=H_R/\langle \gamma \rangle$,  
$A_2=L/\langle \gamma \rangle$ and 
$H_L/\langle \gamma \rangle$. 
Here $\Mod(A_1)=\Mod(A_3)=\pi/2l_X(c)$ and 
$\Mod(A_2)=\pi/l_X(c)$. 
From now on, we regard $A=A_1 \amalg A_2 \amalg A_3$ 
as a subset of $\Sigma^0_c$. 
By definition of grafting, 
it hold that 
$A \cap \Omega^+_{\Sigma^0_c}=A_1 \sqcup A_3$  
and that $\Omega^-_{\Sigma^0_c}=A_2$ on $\Sigma^0_c$, 
where $\Omega^+_{\Sigma^0_c} \subset \Sigma^0_c$ 
(resp. $\Omega^-_{\Sigma^0_c} \subset \Sigma^0_c$)  
is the pull-back 
of the upper half plane $H=\Omega^+_{\Gamma^0}$ 
(resp.  the lower half plane $L=\Omega^-_{\Gamma^0}$).  

Let $\rho \in B_X \cup B_{\bar{X}}$ be a quasiconformal deformation 
of $\rho^0=B(X,\bar{X})$  
which is induced by a quasiconformal map 
$q:\widehat{\mathbb C} \to \widehat{\mathbb C}$.  
If $\rho$ is contained in $B_X$ (resp. $B_{\bar{X}}$)   
we may assume that the support of the Beltrami differential of 
$q$ is contained in 
the lower half plane $L=\Omega^-_{\Gamma^0}$ 
(resp. the upper half plane $H=\Omega^+_{\Gamma^0}$). 
Associated to the quasiconformal deformation $\rho$ of $\rho^0$, 
we have a quasiconformal deformation $\Sigma_c=\Psi_{c}(\rho)$ of 
$\Sigma^0_c=\Psi_{c}(\rho^0)$ 
and a quasiconformal map 
$\varphi:\Sigma^0_c \to \Sigma_c$. 
It is important to notice that the map $\varphi$ takes 
$\Omega^+_{\Sigma^0_c}$ (resp. $\Omega^-_{\Sigma^0_c}$) 
conformally into $\Sigma_c$ if $\rho \in B_X$ 
(resp. $\rho \in B_{\bar{X}}$). 
Here we claim that 
the annulus $\varphi(A) \subset \Sigma_c$ always satisfies 
$\Mod(\varphi(A)) \ge \pi/l_X(c)$. 
In fact, if $\rho \in B_X$ then  
the the map $\varphi$ takes $A_1 \sqcup A_3 \subset \Omega^+_{\Sigma^0_c}$ 
conformally into $\varphi(A) \subset \Sigma_c$ and thus 
$\Mod(\varphi(A)) \ge \Mod(\varphi(A_1))+\Mod(\varphi(A_3))=\pi/l_X(c)$. 
Similarly, if $\rho \in B_{\bar{X}}$ 
then $\varphi$ takes $A_2 \subset \Omega^-_{\Sigma^0_c}$ conformally into 
$\varphi(A) \subset \Sigma_c$ and thus 
$\Mod(\varphi(A)) \ge \Mod(\varphi(A_2))=\pi/l_X(c)$. 
Recall that $X'$ denote the underlying conformal structure of 
$\Sigma_c=\Psi_c(\rho)$. 
By Lemma 4.1, we have 
$\Mod(\varphi(A)) \, l_{X'}(c)^2 \le \Area(X')=4\pi(g-1)$. 
By combining this inequality with the inequality 
$\Mod(\varphi(A)) \ge \pi/l_X(c)$, the result follows 
in the case where $\lambda=c$. 

Next suppose that 
$\lambda \in \mln$ contains a weighted simple closed curve 
$kc$ with $k \in {\mathbb N}$. 
Then $\Sigma_\lambda=\Psi_\lambda(\rho)$ contains 
an annulus $\tilde{A}$ 
which is a union of succeeding 
$k$ parallel annuli each of whose modulus $\ge \pi/l_X(c)$. 
Thus $\Mod(\tilde{A}) \ge k\pi/l_X(c)$ and hence  
we have obtained the desired inequality in general. 
\end{proof}
From Theorem 4.2, we see that the both 
$\pi \circ \Psi_\lambda (B_X)$ and 
$\pi \circ \Psi_\lambda (B_{\bar{X}})$ 
are proper subsets of $T(S)$ for every non-zero $\lambda \in \mln$. 
In fact if $kc \subset \lambda$, they are contained in 
the proper subset  
$$
\{X' \in T(S) : l_{X'}(c) \le \sqrt{4(g-1)\, l_X(c)/k}\}
$$ 
of $T(S)$.  
Compare this fact with a result which was independently obtained by 
Gallo \cite{Ga} and Tanigawa \cite{Ta};  
they showed that 
the map $\pi \circ \Psi_{\lambda}$ takes 
the Fuchsian space 
${\mathcal F}=\{B(X,{\bar X}) \in \QF : X \in T(S)\} \subset \QF$ 
bijectively onto $T(S)$ for every $\lambda \in \mln$. 
In addition, we remark that both 
$\pi \circ \Psi_\lambda (B_X)$ and 
$\pi \circ \Psi_\lambda (B_{\bar{X}})$ 
 are non-precompact in $T(S)$ for every non-zero $\lambda \in \mln$.  
In fact if $c \subset \lambda$, we can make 
the length $l_{X'}(c)$ arbitrarily small 
by varying $\rho=B(X,\bar{Y})$ in $B_X$ 
so that $l_{\bar{Y}}(c) \to 0$ 
or $\rho=B(Y,\bar{X})$ in $B_{\bar{X}}$ so that $l_Y(c) \to 0$.

\subsection{Finiteness theorem}
As a consequence of Theorem 4.2, we obtain the following 
\begin{theorem}[Finiteness]
Suppose that a sequence $\Sigma_n \in Q(S)$  
of quasifuchsian projective structures 
converges to $\Sigma_\infty \in P(S)$ and that
the sequence $\rho_{\Sigma_n} \in \QF$ 
of their holonomies 
converges standardly to 
$\rho_{\Sigma_\infty} \in \partial^{\pm}\QF$. 
Then there is a finite union 
$\Q_{\lambda_1} \sqcup \cdots \sqcup \Q_{\lambda_l}$
of components of $Q(S)$ 
in which $\Sigma_n$ are contained for all $n$. 
\end{theorem}
Instead proving this theorem directly, 
we prove the following theorem, 
which is a contraposition of Theorem 4.3. 
\begin{theorem}
Suppose that a sequence $\rho_n \in \QF$  
converges  standardly to $\rho_\infty \in \partial^{\pm}\QF$ and 
$\{\lambda_n\}$ is a sequence of distinct element of $\mln$. 
Then the sequence 
$\Psi_{\lambda_n}(\rho_n)$ converges to $\infty$  
in $\widehat{P}(S)$ as $n \to \infty$, that is, 
the set 
$\{\Psi_{\lambda_n}(\rho_n):n \in {\mathbb N}\}$ 
has no accumulation point in $P(S)$. 
\end{theorem}
\begin{proof}
Set $\rho_n=B(X_n,\bar{Y}_n)$. 
We first assume that $\rho_\infty \in\partial^+ \QF$. 
Then $\rho_\infty \in \partial B_X$ for some $X \in T(S)$ and 
$X_n \to X$ in $T(S)$.  
One can choose a sequence of weighted simple closed curves 
$k_n c_n \subset \lambda_n$ which satisfies 
$k_n \cdot l_X(c_n) \to \infty$ as $n \to \infty$. 
Let $X_n' \in T(S)$ denote the underlying conformal structure 
of $\Psi_{\lambda_n}(\rho_n)$. 
Then it follows from Theorem 4.2 that 
$$
\frac{l_{X_n'}(c_n)}{l_{X_n}(c_n)} 
\le \sqrt{\frac{4(g-1)}{k_n \cdot l_{X_n}(c_n)}} 
$$
and that the Teichm\"{u}ller distance 
$d_{T(S)}(X_n,X'_n)$ 
between $X_n$ and $X'_n$ tends to $\infty$ as $n \to \infty$. 
Since $X_n \to X$ in $T(S)$, 
$d_{T(S)}(X,X_n')$ also tends to $\infty$,  
which implies that $\Psi_{\lambda_n}(\rho_n) \to \infty$ in $\widehat{P}(S)$.  
The proof for $\rho_\infty \in \partial^- \QF$ is 
completely parallel. 
\end{proof}

\section{Continuity}
The following is our main theorem in this paper. 
\begin{theorem}[Continuity]
Let 
$\rho_n \in \QF$ be a sequence which converges standardly to 
$\rho_\infty \in \partial^{\pm} \QF$. 
Then the sequence $\Psi_{\lambda}(\rho_n)$ converges to 
$\Psi_{\lambda}(\rho_\infty)$ in $\widehat{P}(S)$ 
for every $\lambda \in \mln$. 
\end{theorem}
We devote this section to the proof of the above theorem, 
which is divided into two cases; the case where 
$\Psi_{\lambda}(\rho_\infty) \in P(S)$ and that 
where $\Psi_{\lambda}(\rho_\infty)=\infty$. 
Throughout this proof, 
we let $s$ denote either $+$ or $-$ for which 
$\rho_\infty \in \partial^s \QF$. 
\subsection{Case where $\Psi_{\lambda}(\rho_\infty) \in P(S)$} 
Let $\Phi:U \to V, \, \rho_\infty \mapsto \Psi_\lambda(\rho_\infty)$ 
be a univalent local branch of $\hol^{-1}$  
from a neighborhood $U$ of $\rho_\infty$ 
to that $V$ of $\Psi_\lambda(\rho_\infty)$ and set 
$\Sigma_n:=\Phi(\rho_n)$ and 
$\Sigma_\infty:=\Phi(\rho_\infty)=\Psi_\lambda(\rho_\infty)$. 
Then the sequence $\Sigma_n$ 
converges to $\Sigma_\infty$ in $P(S)$. 
We claim in 
Lemma 5.2 below that $\Sigma_n \in \Q_\lambda$ for all large enough $n$. 
Assuming this lemma, we have 
$\Sigma_n=\Psi_{\lambda}(\rho_n)$ 
and thus obtain the desired convergence 
$\Psi_{\lambda}(\rho_n)=\Sigma_n \to \Sigma_\infty=\Psi_{\lambda}(\rho_\infty)$.  
\begin{lemma}
$\Sigma_n \in \Q_\lambda$ for all large enough $n$. 
\end{lemma} 

The rest of this subsection is devoted to the proof of Lemma 5.2. 
The idea of the proof is as follows: 
by Lemma 3.6, we have already known the topological type of the set 
$\Omega^s_{\Sigma_\infty} \subset \Sigma_\infty=\Psi_\lambda(\rho_\infty)$ 
in relation to $\lambda$.  
As we will see in Lemma 5.5, the set 
$\Omega^s_{\Sigma_\infty} \subset \Sigma_\infty$ coincides with 
the essential part of the Carath\'{e}odory limit of 
the sequence $\Omega^s_{\Sigma_n} \subset \Sigma_n$ as $n \to \infty$. 
It is essential that the sequence 
$\Sigma_n$ is contained in a finite union of 
components of $Q(S)$ by Theorem 4.3. 
Then we will show that for any subsequence 
$\{\Sigma_{n_j}\}_{j=1}^\infty$ of $\{\Sigma_n\}_{n=1}^\infty$ 
contained in a component of $Q(S)$, 
the topological types of $\Omega^s_{\Sigma_{n_j}} \subset \Sigma_{n_j}$ 
is the same to that of $\Omega^s_{\Sigma_\infty} \subset \Sigma_\infty$ 
for all $j$,  
which implies that $\Sigma_n \in \Q_\lambda$ for all $n$. 

Now we fill in the details.  
We may assume that the sequence 
$\Gamma_n=\rho_n(\pi_1(S))$ 
of quasifuchsian groups converges geometrically to 
some Kleinian group $\widehat{\Gamma}$, which 
contains the algebraic limit 
$\Gamma_\infty=\rho_\infty(\pi_1(S))$. 
Then Kerckhoff and Thurston \cite[Corollary 2.2]{KeTh} showed 
that the sequence $\Omega_{\Gamma_n}$ 
of the regions of discontinuity of $\Gamma_n$ 
converges to the region of discontinuity 
$\Omega_{\widehat{\Gamma}}$ of $\widehat{\Gamma}$ 
in the sense of Carath\'{e}odory. 
Here $\Omega_{\Gamma_n} \to \Omega_{\widehat{\Gamma}}$ in this sense 
if (i) for every compact $K \subset \Omega_{\widehat{\Gamma}}$, we have 
$K \subset \Omega_{\Gamma_n}$ for all large enough $n$ 
and if (ii) for any open set 
$U \subset \Omega_{\Gamma_{n_j}} \, (n_j \to \infty)$, we have 
$U \subset \Omega_{\widehat{\Gamma}}$. 
Moreover $\Omega_{\widehat{\Gamma}}$ decomposes  
into two $\widehat{\Gamma}$-invariant regions  
$\Omega^+_{\widehat{\Gamma}}$, $\Omega^-_{\widehat{\Gamma}}$, 
each of which are the Carath\'{e}odory limits of 
$\Omega^+_{\Gamma_n}$, $\Omega^-_{\Gamma_n}$, respectively. 
Since the convergence $\rho_n \to \rho_\infty \in \partial^s \QF$ 
is standard, 
the argument in the proof of \cite[Proposition 2.3]{KeTh} 
reveals that 
$\Omega^s_{\Gamma_\infty}$ 
is a component of $\Omega^s_{\widehat{\Gamma}}$, that 
$\Omega^s_{\widehat{\Gamma}} 
=\bigcup_{\delta \in \widehat{\Gamma}} \delta(\Omega^s_{\Gamma_\infty})$, 
and that $\delta \in \widehat{\Gamma}$ fixes 
$\Omega^s_{\Gamma_\infty}$ if and only if $\delta \in \Gamma_\infty$. 
We put them together in the following 
\begin{lemma}[Kerckhoff-Thurston]
We have 
$\Omega^s_{\widehat{\Gamma}} 
=\bigsqcup_{\delta \in \widehat{\Gamma}/\Gamma_\infty} 
\delta(\Omega^s_{\Gamma_\infty})$. 
\end{lemma}
Recall that  
$\Omega^{s}_{\Sigma_n}$ denotes 
the pull-back of $\Omega^s_{\Gamma_n}$ in $\Sigma_n$. 
Similarly $\Omega^{s}_{\Sigma_\infty}$ is  
the pull-back of $\Omega^s_{\Gamma_\infty}$ in  $\Sigma_\infty$. 
In addition, we let 
$ \widehat{\Omega}^{s}_{\Sigma_\infty}$ denotes the pull-back of   
$\Omega^{s}_{\widehat{\Gamma}}$ in $\Sigma_\infty$.  
By the same argument of \cite[Lemma 3.3]{It1}, we obtain the following 
\begin{lemma}
The sequence 
$\Omega^s_{\Sigma_n} \subset \Sigma_n$ converges to $
\widehat{\Omega}^s_{\Sigma_\infty} \subset \Sigma_\infty$ 
in the sense of Carath\'{e}odory.   
More precisely, 
there exist $K_n$-bi-Lipshiz diffeomorphisms 
$\psi_n:\Sigma_\infty \to \Sigma_n$ 
between hyperbolic surfaces 
(consistent with their markings) such that $K_n \to 1$ as $n \to \infty$ 
and that the sequence 
$\psi_n^{-1}(\Omega^s_{\Sigma_n}) \subset \Sigma_\infty$ 
converges to 
$\widehat{\Omega}^s_{\Sigma_\infty} \subset \Sigma_\infty$ 
in the sense of Carath\'{e}odory. 
\end{lemma} 
We claim in the following lemma that 
$\Omega^s_{\Sigma_\infty} \subset \Sigma_\infty$ carries the full 
fundamental group of 
$\widehat{\Omega}^s_{\Sigma_\infty} \subset \Sigma_\infty$.  
\begin{lemma}
Each connected component of 
$\widehat{\Omega}^s_{\Sigma_\infty} -\Omega^s_{\Sigma_\infty}$ 
is homotopically trivial in $\Sigma_\infty$. 
\end{lemma}
\begin{proof}
Suppose for contradiction that 
there exists a connected component $\omega$ of 
$\widehat{\Omega}^s_{\Sigma_\infty} -\Omega^s_{\Sigma_\infty}$ 
which contains an essential closed curve $\alpha$.  
Let $\tilde{\omega} \subset \widetilde{\Sigma}_\infty$ be 
a lift (i.e., a connected component of the inverse image) 
of $\omega \subset \Sigma_\infty$ and 
$\tilde{\alpha} \subset \widetilde{\Sigma}_\infty$ 
a lift of the curve $\alpha \subset \omega$ 
which is contained in $\tilde{\omega}$.   
Let $\Omega$ be the connected component of 
$\Omega^s_{\widehat{\Gamma}}-\Omega^s_{\Gamma_\infty}$ 
which contains the developing image 
$f_{\Sigma_\infty}(\tilde \omega) \subset \widehat{\mathbb C}$. 
Let $\gamma=\rho_\infty(\alpha)$ be an element of $\Gamma_\infty$ 
which fixes $f_{\Sigma_\infty}(\tilde{\alpha})$, and hence $\Omega$.  

Now we see that $\gamma$ is loxodromic 
from Definition 3.3 and Lemma 3.5. 
In fact if $\rho_\infty \in \partial^+ \QF$, then 
each connected component of 
$\Sigma_\infty-\Omega^+_{\Sigma_\infty}$ is an annulus 
whose core curve is homotopic to a component of $\lambda$. 
Hence the curve $\alpha$ in 
$\omega \subset \Sigma_\infty-\Omega^+_{\Sigma_\infty}$ is homotopic 
into $\lambda$ and hence 
is not homotopic into $\para(\rho_\infty)$. 
On the other hand, 
if $\rho_\infty \in \partial^- \QF$ then each component of 
$\para(\rho_\infty)$ must traverse $\lambda$, 
and thus $\Omega^-_{\Sigma_\infty}$, in $\Sigma_\infty$. 
Therefore the curve $\alpha$ in $\Sigma_\infty-\Omega^-_{\Sigma_\infty}$ 
is not homotopic into $\para(\rho_\infty)$. 
 
From Lemma 5.3, there exists 
$\delta \in \widehat{\Gamma}-\Gamma_\infty$ such that 
$\Omega=\delta(\Omega^s_{\Gamma_\infty})$. 
Since $\gamma$ fixes the component $\Omega$,  
$\delta^{-1} \gamma \delta$ fixes the component 
$\Omega^s_{\Gamma_\infty}$,  which implies that 
$\delta^{-1} \gamma \delta \in \Gamma_\infty$.  
But this contradicts to the following lemma 
(Lemma 2.4 in \cite{ACCS}). 
\end{proof}
\begin{lemma}[Anderson-Canary-Culler-Shalen]
Let $G$ be a finitely generated group. 
Let $\rho_n:G \to \Gamma_n \subset \psl$ be a sequence of 
discrete faithful representations converging algebraically to 
$\rho_\infty:G \to \Gamma_\infty$. 
Moreover, assume that $\Gamma_n$ converges to 
$\widehat{\Gamma}$ geometrically. 
Then for any $\delta \in \widehat{\Gamma}-\Gamma_\infty$,  
$\Gamma_\infty \cap \delta^{-1} \Gamma_\infty \delta$ 
is $\{\id\}$ or rank-one parabolic subgroup.  
\end{lemma}
Let us combine the above observations 
to complete the proof of Lemma 5.2, which states that 
$\Sigma_n \in {\mathcal Q}_{\lambda}$ for all large enough $n$. 
Note that the sequence 
$\Sigma_n$ is contained in a finite union of 
components of $Q(S)$ by Theorem 4.3. 
Therefore to show Lemma 5.2, 
we only need to show that $\lambda'=\lambda$ by 
assuming that the sequence $\Sigma_n$ is contained 
in a component $\Q_{\lambda'}$ of $Q(S)$. 

Let $\rho_0 \in \QF$ be a fixed element 
and set $\Sigma_0=\Psi_{\lambda'}(\rho_0)$.
Then every $\Sigma_n$ is a quasiconformal deformation of $\Sigma_0$ 
in $\Q_{\lambda'}$, 
which is induced by the quasiconformal deformation $\rho_n$ 
of $\rho_0$, and which induces  
the quasiconformal map  
$\varphi_n: \Sigma_0 \to \Sigma_n$. 
Note that the map $\varphi_n$ is consistent with their markings 
and takes $\Omega^s_{\Sigma_0}$ onto $\Omega^s_{\Sigma_n}$. 
Since the convergence $\rho_n \to \rho_\infty \in \partial^s \QF$ 
is standard, the dilatation of $\varphi_n|_{\Omega^s_{\Sigma_0}}$ 
tends to $1$, and hence we assume that 
$\varphi_n|_{\Omega^s_{\Sigma_0}}:\Omega^s_{\Sigma_0} \to \Omega^s_{\Sigma_n}$ 
is conformal for every $n$ for simplicity. 

Now let ${\mathcal C}(\Omega^s_{\Sigma_0})$, 
${\mathcal C}(\Omega^s_{\Sigma_\infty})$ denote  
the set of connected components of $\Omega^s_{\Sigma_0}$, 
$\Omega^s_{\Sigma_\infty}$, respectively. 
Then Lemma 3.6 tells us that $\lambda=\lambda'$ if and only if 
there exists a bijection 
$$
g:{\mathcal C}(\Omega^s_{\Sigma_0}) 
\to {\mathcal C}(\Omega^s_{\Sigma_\infty})
$$
which satisfies 
$\pi_1(\omega)=\pi_1(g(\omega))$ for every 
$\omega \in {\mathcal C}(\Omega^s_{\Sigma_0})$. 
Here we identify $\pi_1(\Sigma_0)$ and $\pi_1(\Sigma_\infty)$ 
with $\pi_1(S)$ via their markings. 
Note that 
$\psi^{-1}_n \circ \varphi_n(\Omega^s_{\Sigma_0}) \subset \Sigma_\infty$ 
converge to $\widehat{\Omega}^s_{\Sigma_\infty} \subset \Sigma_\infty$ 
is the sense of Lemma 5.4 and that 
$\pi_1(\omega)=
\pi_1(\psi^{-1}_n \circ \varphi_n(\omega))$ 
hold for all $\omega \in {\mathcal C}(\Omega^s_{\Sigma_0})$ and $n$. 
We first define a map 
$h:{\mathcal C}(\Omega^s_{\Sigma_\infty}) \to {\mathcal C}(\Omega^s_{\Sigma_0})$, 
which will turn out to be the inverse of $g$.  
Let $\omega' \in {\mathcal C}(\Omega^s_{\Sigma_\infty})$ and let 
$K \subset \omega'$ be a compact subset  
such that $\pi_1(K)=\pi_1(\omega')$.  
Then there exists $\omega \in {\mathcal C}(\Omega^s_{\Sigma_0})$ 
such that $K \subset \psi^{-1}_n \circ \varphi_n(\omega)$ 
for all large enough $n$,  
and hence that  
$\pi_1(K)=\pi_1(\omega') \subset \pi_1(\omega)$. 
Now we let $h(\omega')=\omega$. 
Next we define the map 
$g:{\mathcal C}(\Omega^s_{\Sigma_0}) \to 
{\mathcal C}(\Omega^s_{\Sigma_\infty})$ 
as follows: 
let $\omega \in {\mathcal C}(\Omega^s_{\Sigma_0})$. 
By passing to a subsequence if necessary, 
the sequence 
$\psi^{-1}_n \circ \varphi_n(\omega) \subset \Sigma_\infty$ 
converges  to an open subset $\hat{\omega}$ of 
$\widehat{\Omega}_{\Sigma_\infty}^s \subset \Sigma_\infty$. 
We claim that $\hat{\omega}$ contains 
exactly one essential component $\omega'$ of 
$\widehat{\Omega}^s_{\Sigma_\infty} \subset \Sigma_\infty$,   
which satisfies $\pi_1(\omega) \subset \pi_1(\omega')$. 
We postpone the proof of this claim and continue the argument. 
We have $\omega' \in {\mathcal C}(\Omega^s_{\Sigma_\infty})$ 
by Lemma 5.5 and let $g(\omega)=\omega'$. 
Then it follows immediately that 
$g$ is a bijective map with $g^{-1}=h$ and hence that $\lambda=\lambda'$.  
We have completed the proof of Lemma 5.2, and 
hence that of Theorem 5.1 in the case where 
$\Psi_\lambda(\rho_\infty) \in P(S)$. 

We only need to show the above claim. 
Suppose for contradiction that 
$\hat{\omega} \subset \Sigma_\infty$ is trivial 
or contains more that two essential connected components.  
Since $\psi^{-1}_n \circ \varphi_n(\omega)$ converge to $\hat{\omega}$, 
$\varphi_n(\omega) \subset \Sigma_n$ then 
become more and more constricted, 
that is, there exists an arc 
$\alpha_n \subset \varphi_n(\omega)$ joining 
two distinct components of $\varphi_n(\partial \omega)$ 
whose hyperbolic length tends to $0$ as $n \to \infty$. 
But since the map 
$\varphi_n|_{\Omega^s_{\Sigma_0}}:\Omega^s_{\Sigma_0} \to 
\varphi_n(\Omega^s_{\Sigma_0})$ 
is conformal, 
this contradicts to Lemma 5.7 below. 
Thus we have proved the claim. 
\begin{lemma}[{\cite[Lemma 4.4]{It1}}]
Let $R$ be a hyperbolic surface, $c \subset R$ a simple closed curve  
and $A \subset R$ is an annulus whose core curve is homotopic to $c$. 
Then there is a positive constant $C>0$,  
depending only on $l_R(c)$ and $\Mod(A)$,  
such that $l_R(\alpha) \ge C$ holds  
for every arcs $\alpha \subset A$ 
joining points in distinct components of $\partial A$. 
\end{lemma}

\subsection{Case where $\Psi_{\lambda}(\rho_\infty) = \infty$} 
We first prepare the following 
\begin{lemma}
Let $\lambda \in \mln$ and 
let $\rho_n=B(X_n,\bar{Y}_n) \in \QF$ be a sequence 
which converges standardly to $\rho_\infty \in \partial^\pm \QF$ 
and which satisfies the property as in Corollary 2.4. 
Suppose that $\Psi_{\lambda}(\rho_\infty) = \infty$. 
Then the sequence $\Psi_{\lambda}(\rho_n) \to \infty$ 
in $\widehat{P}(S)$. 
\end{lemma}
\begin{proof}
We first suppose that $\rho_\infty \in \partial^+ \QF$. 
Since $(\lambda,\rho_\infty)$ is not admissible, 
$\lambda$ and $\para(\rho_\infty)$ 
have parallel components in common 
and let $c$ be one of such components.  
Since $l_{\bar{Y}_n}(c) \to 0$, 
there are annular neighborhoods $A_n$ 
of $c$ in $\bar{Y}_n$ such that $\Mod(A_n) \to \infty$
as $n \to \infty$. 
We set $\Sigma_n=\Psi_{\lambda}(\rho_n)$. 
Then observe that the annulus $A_n$ is 
conformally embedded in 
$\Omega^-_{\Sigma_n} \subset \Sigma_n$ for every $n$. 
Hence the underlying conformal structures $\pi(\Psi_{\lambda}(\rho_n))$ 
of $\Psi_{\lambda}(\rho_n)$ 
diverge in $T(S)$, and hence $\Psi_{\lambda}(\rho_n) \to \infty$ 
in $\widehat{P}(S)$. 

Next we suppose that $\rho_\infty \in \partial^- \QF$. 
Then there is a component $c$ of $\para(\rho_\infty)$ 
which do not intersect $\lambda$ essentially. 
Since $l_{X_n}(c) \to 0$, 
there are annular neighborhoods $A_n$ 
of $c$ in $X_n$ such that $\Mod(A_n) \to \infty$ 
as $n \to \infty$. 
Then observe that the annulus $A_n$ is 
conformally embedded in 
$\Omega^+_{\Sigma_n} \subset \Sigma_n=\Psi_{\lambda}(\rho_n)$ for every $n$. 
By the same argument as above, we see that 
$\Psi_{\lambda}(\rho_n) \to \infty$ 
in $\widehat{P}(S)$. 
\end{proof}
Now we back to the proof of Theorem 5.1 in the case 
where $\Psi_\lambda(\rho_\infty)=\infty$; 
we intend to show that 
$\Psi_{\lambda}(\rho_n) \to \Psi_\lambda(\rho_\infty)=\infty$ 
for every standard convergent sequence 
$\QF \ni \rho_n \to \rho_\infty \in \partial^\pm \QF$. 
Set $\Sigma_n=\Psi_\lambda(\rho_n)$.  
To obtain a contradiction, 
we suppose that the sequence $\Sigma_n$ has a subsequence 
(which is denoted by the same symbols) 
converging to some $\Sigma_\infty \in P(S)$.  
In addition, let 
$\rho_n' \in \QF$ be a sequence 
converging to $\rho_\infty$ as in Lemma 5.8 and 
let $\Phi:U \to P(S), \ \rho_\infty \mapsto \Sigma_\infty$ 
be a univalent local branch of $\hol^{-1}$.  
Then we have $\Sigma_n=\Phi(\rho_n)$ 
for all large enough $n$. 
Note that both 
$\Sigma_n=\Phi(\rho_n)$ 
and $\Sigma_n':=\Phi(\rho_n')$ converge to 
$\Sigma_\infty$.  
From Theorem 4.3, 
the sequence $\Sigma_n'$ has a subsequence 
(also denoted by the same symbols) 
which is contained in a component ${\mathcal Q}_{\mu}$ of $Q(S)$ 
for some $\mu \in \mln$. 
Then we have 
$\Sigma_n'=\Psi_{\mu}(\rho_n')$, and thus 
$\Psi_{\mu}(\rho_n')=\Sigma_n' \to \Sigma_\infty \ne \infty$, 
which implies that 
$\Psi_{\mu}(\rho_\infty) \ne \infty$ by Lemma 5.8 and that 
$\lambda \ne \mu$. 
Since both $\rho_n,\,\rho'_n$ converge standardly to $\rho_\infty$, 
Theorem 5.1 for admissible pairs guarantees that 
both $\Psi_\mu(\rho_n),\,\Psi_\mu(\rho'_n)$ 
convege to $\Psi_\mu(\rho_\infty)$. 
On the other hand, since $\Psi_{\mu}(\rho_n')=\Sigma_n' \to \Sigma_\infty$, 
we have $\Psi_\mu(\rho_n) \to \Sigma_\infty$. 
Now we obtain two sequences 
$\Psi_\lambda(\rho_n), \, \Psi_\mu(\rho_n)$ 
both of which converge to $\Sigma_\infty$. 
But since $\lambda \ne \mu$, this contradicts to the fact that 
the map $\hol$ is a local homeomorphism. 
Now we have completed the proof of Theorem 5.1.

\section{Applications}
We collect in this section some results 
which are obtained as consequences of Theorem 5.1. 
\subsection{Grafting theorem for boundary $b$-groups}
The following statement is conjectured by Bromberg in \cite{Br}. 
\begin{theorem}
For every boundary $b$-groups $\rho \in  \partial^\pm \QF$, 
all projective structure with holonomy $\rho$ 
are obtained by grafting of $\rho$, that is, we have 
$$
\hol^{-1}(\rho)=\{\Psi_{\lambda}(\rho) : 
\lambda \in \mln, \Psi_{\lambda}(\rho) \in P(S)\}.
$$
\end{theorem}
\begin{proof} 
Let $\rho \in \partial^\pm \QF$ 
and $\Sigma \in \hol^{-1}(\rho)$.  
We only need to show that 
$\Sigma=\Psi_{\lambda}(\rho)$ for some 
$\lambda \in \mln$.  
Let $\rho_n \in \QF$ be a sequence   
converging standardly to $\rho$.  
Since the map $\hol$ is a local homeomorphism, 
there exists a convergent sequence $\Sigma_n \to \Sigma$ 
of projective structures with $\hol(\Sigma_n)=\rho_n$ for all $n$.  
From Theorem 4.3, 
there exists a subsequence $\{\Sigma_{n_j}\}_{j=1}^\infty$ 
of the sequence $\{\Sigma_n\}_{n=1}^\infty$ 
which is contained in a component 
$\Q_{\lambda}$ of $Q(S)$. 
Then we have $\Sigma_{n_j}=\Psi_{\lambda}(\rho_{n_j})$ for all $j$ 
and thus 
$\Psi_{\lambda}(\rho_{n_j})=\Sigma_{n_j} \to \Sigma$ as $j \to \infty$. 
On the other hand, 
since the convergence $\rho_n \to \rho$ is standard, 
Theorem 5.1 implies that 
$\Psi_\lambda(\rho_n) \to \Psi_\lambda(\rho)$ in $\widehat{P}(S)$.  
Hence
$\Sigma=\Psi_{\lambda}(\rho)$ 
and the result follows. 
\end{proof}
From Bromberg's result (Theorem 2.2), we obtain the following 
\begin{corollary}
Suppose that $\rho \in AH(S)$ is a $b$-group 
with no parabolics. 
Then we have $\hol^{-1}(\rho)=\{\Psi_{\lambda}(\rho) : \lambda \in \mln\}$. 
\end{corollary}
\subsection{Analytic continuations and lifts of paths}
\begin{theorem}
Suppose that $\Psi_\lambda(\rho) \in P(S)$ for 
$\lambda \in \mln$ and $\rho \in \partial^\pm \QF$. 
Then there exists a path 
$\alpha:[0,1] \to R(S)$ with 
$\alpha(0) \in \QF$ and $\alpha(1)=\Psi_\lambda(\rho)$  
along which the grafting map 
$\Psi_\lambda:\QF \to P(S)$ 
is continued analytically 
to a univalent local branch 
$\Phi:U \to P(S), \, \rho \mapsto \Psi_{\lambda}(\rho)$ 
of $\hol^{-1}$ which is 
defined on some neighborhood $U$ of $\rho$. 
\end{theorem}
\begin{proof}
We may assume that $\rho \in \partial^+ \QF$, and hence that 
$\rho \in \partial B_X$ for some $X \in T(S)$.  
Let $D \subset T(S)$ de a bounded domain with $X \in D$ 
and set $B_D =\sqcup_{X' \in D}B_{X'}$.
We claim that there exists a neighborhood 
$U$ of $\rho$ such that 
$\Phi \equiv \Psi_\lambda$ holds on $U \cap B_D$.  
For otherwise there is a 
sequence $\rho_n \in B_D$ converging to $\rho$ for which 
$\Phi(\rho_n) \ne \Psi_\lambda(\rho_n)$ hold for all $n$. 
On the other hand, Theorem 5.1 implies that 
both sequences $\Phi(\rho_n), \Psi_\lambda(\rho_n)$ 
converge to $\Phi(\rho)=\Psi_\lambda(\rho)$, 
which contradicts to the fact that $\Phi$ is univalent.   
Now let  
$\alpha:[0,1] \to R(S)$ be an arc such that 
$\alpha(0) \in \QF$, $\alpha(1)=\rho$ and that 
$\alpha([0,1]) \subset B_D \cup U$. 
Then $\Psi_\lambda:\QF \to P(S)$ is continued 
analytically to $\Phi:U \to P(S)$ along $\alpha$. 
\end{proof}
Let $\alpha:[0,1] \to R(S)$ be a path as in Theorem 6.3. 
Then the above theorem implies that there exists the lift 
$\tilde{\alpha}:[0,1] \to P(S)$ of $\alpha$ 
via the holonomy map  
with the starting point $\Psi_\lambda(\alpha(0))$. 
On the other hand, 
we have the following 
\begin{theorem}
Suppose that $\Psi_\lambda(\rho)=\infty$ for $\lambda \in \mln$ 
and $\rho \in \partial^\pm \QF$. 
Then there exists a path $\alpha:[0,1] \to R(S)$ 
of which there is no lift to $P(S)$ 
with the starting point $\Psi_\lambda(\alpha(0))$. 
\end{theorem}
\begin{proof}
We may assume that $\rho \in \partial^+ \QF$, and hence that 
$\rho \in \partial B_X$ for some $X \in T(S)$.  
Let take a path $\alpha:[0,1] \to R(S)$ 
with $\alpha(0) \in \QF$ and $\alpha(1)=\rho$ 
for which there exists a sequence 
$0<t_n<1, \, t_n \to 1$ such that 
$\alpha(t_n) \in B_D =\sqcup_{X' \in D}B_{X'}$ 
for some bounded domain $D \subset T(S)$ with $X \in D$.  
Then the sequence $\alpha(t_n) \in B_D$ 
converges standardly to $\rho$, and 
thus $\Psi_{\lambda}(\alpha(t_n)) \to \infty$ by Theorem 5.1. 
Hence there is no lift of $\alpha$ to $P(S)$ 
with the starting point $\Psi_\lambda(\alpha(0))$. 
\end{proof}

\subsection{Obstructions for $\hol$ to be a covering map}
We first fix our terminology. 
\begin{definition}
Let ${\mathcal X}, {\mathcal Y}$ be topological spaces.  
A continuous map  
$f:{\mathcal Y} \to {\mathcal X}$ 
is said to be a {\it covering map} or a {\it weak-covering map} if 
for any $x \in {\mathcal X}$ there is a neighborhood $U$ of $x$ such that 
$f|_V:V \to U$ is a homeomorphism 
\begin{itemize}
  \item
  for any connected component $V$ of $f^{-1}(U)$,  or 
  \item  
for any $y \in f^{-1}(x)$ and some neighborhood $V$ of $y$,  
\end{itemize}
respectively. 
\end{definition}
Although the holonomy map $\hol:P(S) \to R(S)$ 
is a local homeomorphism and 
the map $\hol|_{Q(S)}:Q(S) \to \QF$ is a covering map, 
Hejhal showed the following 
\begin{theorem}[{\cite[Theorem 8]{He}}]
The holonomy map $\hol:P(S) \to R(S)$ is not a covering map 
onto its image. 
\end{theorem}
We explain this fact in our context. 
Let $\rho \in \partial^\pm \QF$ with 
$\para(\rho) \ne \emptyset$. 
Then there are $\lambda,\mu \in \mln$ 
such that $\Psi_\lambda(\rho) \in P(S)$ 
and $\Psi_\mu(\rho)=\infty$.  
One can find a path $\alpha:[0,1] \to R(S)$ 
with $\alpha(0) \in \QF$ and $\alpha(1)=\rho$ 
which satisfies both conditions in the proofs 
of Theorems 6.3 and 6.4. 
Then there is the lift of the path $\alpha$ 
with the starting point $\Psi_\lambda(\alpha(0))$ 
but there is no lift with the starting point $\Psi_\mu(\alpha(0))$, 
which implies that the map $\hol$ is not a covering map. 
We remark that the idea used in the argument above is the same 
to that of Hejhal \cite{He}, 
but he made use of a path $\alpha:[0,1] \to R(S)$ of 
Shottky representations.  

Moreover we claim in Corollary 6.8 below that the covering map 
$\hol|_{Q(S)}:Q(S) \to \QF$ is not 
extended any more even in the sense of weak-covering. 
The essential observation is the following 
\begin{theorem}
Let $\rho \in \partial^\pm \QF$. 
For every $\lambda \in \mln$ with $\Psi_\lambda(\rho) \in P(S)$, 
suppose that 
$\Phi_\lambda:U_\lambda \to P(S), \, \rho \mapsto \Psi_\lambda(\rho)$ 
is a univalent local branch of $\hol^{-1}$ 
defined on a neighborhood $U_\lambda$ of $\rho$. 
Then $\rho$ can not be an interior point 
of $\bigcap_\lambda U_\lambda$, 
where the intersection is taken over all 
$\lambda \in \mln$ with $\Psi_\lambda(\rho) \in P(S)$. 
\end{theorem}
\begin{proof}
We may assume that $\rho \in \partial B_X$ for some $X \in T(S)$.  
Then there exists a sequence $\rho_n$ 
of maximal cusps converging to $\rho$ in $\partial B_X$ 
for which the parabolic locus $\lambda_n=\para(\rho_n) \in \mln$ 
has no parallel component in common 
with $\para(\rho)$ (which is possibly empty) for every $n$ 
(see \cite[Theorem 5.5]{It2}). 
Since $\Psi_{\lambda_n}(\rho) \ne \infty$,  
there exists a univalent local branch 
$\Phi_{\lambda_n}:U_{\lambda_n} \to P(S),
\,\rho \mapsto \Psi_{\lambda_n}(\rho)$ 
for every $n$.  
On the other hand, since $\Psi_{\lambda_n}(\rho_n)=\infty$, 
$\rho_n \not\in U_{\lambda_n}$ for every $n$. 
Therefore we have obtained the result. 
\end{proof}

\begin{corollary}
For any open neighborhood $O$ of $\QF$ in $R(S)$, 
the map $\hol|_{\hol^{-1}(O)}:\hol^{-1}(O) \to O$ 
is not a weak-covering map. 
\end{corollary}
\begin{proof}
For any $\rho \in \partial \QF$ and its neighborhood $U$ contained in $O$, 
one can find a boundary 
$b$-group $\rho' \in \partial^\pm \QF$ in $U$. 
Then the result follows from the above theorem.  
\end{proof}

\subsection{Further continuity} 
We extend our definition of standard convergence 
to sequences in $\overline{\QF}$. 
Let $K$ be a compact subset of $T(S)$ and 
$\bar{K} \subset T(\bar{S})$ its complex conjugation. 
We set $B_K=\bigsqcup_{X \in K}B_X$ and 
$B_{\bar{K}}=\bigsqcup_{\bar{X} \in \bar{K}}B_{\bar{X}}$. 
Recall that a sequence $\rho_n \in \QF$ converges 
standardly to $\rho_\infty \in \partial \QF$ 
if there is a compact subset $K \subset T(S)$ 
such that $\rho_n \in B_K \cup B_{\bar{K}}$ for all $n$. 
Similarly we say that 
a sequence $\rho_n \in \overline{\QF}$ converges 
{\it standardly} to $\rho_\infty \in \partial \QF$ 
if there is a compact subset $K \subset T(S)$ 
such that all $\rho_n$ are contained in the 
{\it closure} of $B_K \cup B_{\bar{K}}$ in $R(S)$.  
Note that the limit of a standard convergent sequence 
in $\overline{\QF}$ is contained in $\partial^\pm \QF$ 
and that if a sequence $\rho_n \in \partial \QF$ converges standardly to 
$\partial^s \QF$ then $\rho_n \in \partial^s \QF$ for all large enough $n$, 
where $s$ denotes $+$ or $-$. 
\begin{theorem}
Let 
$\rho_n \to \rho_\infty$ be a standard convergent sequence 
in $\overline{\QF}$. 
Then the sequence $\Psi_{\lambda}(\rho_n)$ converges to 
$\Psi_{\lambda}(\rho_\infty)$ in $\widehat{P}(S)$ for every $\lambda \in \mln$. 
\end{theorem} 
\begin{proof}
We only have to consider the case where $\rho_n \in \partial \QF$ for all $n$. 
Then $\rho_n \to \rho_\infty$ in $\partial^\pm \QF$. 
If $\Psi_\lambda(\rho_\infty) \in P(S)$ then the result 
follows from Theorem 6.3.  
Hence we suppose that $\Psi_\lambda(\rho_\infty)=\infty$. 
(We remark that the argument below also works in the case where 
$\Psi_\lambda(\rho_\infty) \in P(S)$.) 
To obtain a contradiction, 
we assume that 
the sequence $\Psi_\lambda(\rho_n)$ converges to some 
$\Sigma' \in P(S)$. 
Then for each $n$ 
there exists a sequence 
$\rho_{n,k} \in \QF$ converging standardly to 
$\rho_n \in \partial^\pm \QF$ as $k \to \infty$, 
which is taken by $\Psi_\lambda$ to a 
convergent sequence 
$\Psi_\lambda(\rho_{n,k}) \to \Psi_\lambda(\rho_n)$ from Theorem 5.1. 
Then by a diagonal argument, we can choose a sequence 
$\rho_n'$ in $\QF$ which converges standardly to $\rho_\infty$ 
and which satisfies $\Psi_\lambda(\rho_n') \to \Sigma'$. 
On the other hand, 
$\Psi_\lambda(\rho_n') \to \Psi_\lambda(\rho_\infty)=\infty$ by Theorem 5.1, 
which is a contradiction.   
Now we have obtained the result.  
\end{proof}
Let $\overline{\Q_0}$ denote the closure of $\Q_0$ in $P(S)$, 
{\it not} in $\widehat{P}(S)$. 
From the observation in \S 3.3, we see that 
a sequence $\Sigma_n \in \overline{\Q_0}$ converges to 
$\Sigma_\infty \in \partial \Q_0$ in $P(S)$ 
if and only if the sequence $\rho_{\Sigma_n} \in \QF \sqcup \partial^+ \QF$ 
converges standardly to $\rho_{\Sigma_\infty} \in \partial^+ \QF$. 
Hence, as a direct consequence of the above theorem, we obtain the 
following 
\begin{theorem}
The map $\Gr_{\lambda}:\Q_0 \to \Q_{\lambda}$ 
is extended continuously 
to $\Gr_\lambda:\overline{\Q_0} \to \widehat{P}(S)$ 
for each $\lambda \in \mln$,  
where $\overline{\Q_0}$ is the closure of $\Q_0$ in $P(S)$. 
\end{theorem}

\subsection{Discreteness of inverse images}
\begin{theorem}
Let $K$ be a given compact subset of $T(S)$ 
and set $\QF_K=B_K \cup B_{\bar{K}}$. 
Then the inverse image 
$$
\hol^{-1}(\QF_K)=
\bigsqcup_{\lambda \in \mln} \Psi_\lambda (\QF_K)
$$ 
of $\QF_K$ in $P(S)$ is discrete; that is, 
every connected component of $\hol^{-1}(\QF_K)$ 
has open neighborhood which is disjoint from any other one. 
\end{theorem}
\begin{proof}
Suppose for contradiction that there exist  
$\rho \in \overline{\QF_K} \cap \partial^\pm\QF$ 
and $\lambda \in \mln$ such that 
$\Psi_{\lambda}(\rho)$ is a limit of a sequence $\Sigma_n$ in 
$\hol^{-1}(\QF_K)-\Psi_\lambda(\QF_K)$. 
By Theorem 4.3, we may assume that 
the sequence $\Sigma_n$ is contained in a component $\Q_\mu$ of $Q(S)$ 
with $\lambda \ne \mu$. 
Since the sequence $\rho_n=\hol(\Sigma_n) \in \QF_K$ 
converges standardly to $\rho$, 
we have $\Sigma_n =\Psi_\mu(\rho_n) \to \Psi_\mu(\rho)$ by Theorem 5.1. 
But this contradicts the fact that 
$\Psi_\lambda(\rho) \ne \Psi_\mu(\rho)$ (Corollary 3.7).  
\end{proof}
Compare the above theorem with the following 
\begin{theorem}[\cite{It3}]
Any pair of connected components of the inverse image
$$
\hol^{-1}(\QF)=\bigsqcup_{\lambda \in \mln} \Psi_\lambda(\QF) 
$$ 
of $\QF$ 
have intersecting closures in $P(S)$. 
\end{theorem}
\bigskip
\begin{flushleft}
Graduate School of Mathematics, \\
Nagoya University, \\
Nagoya 464-8602, Japan  \\
\texttt{itoken@math.nagoya-u.ac.jp} 
\end{flushleft}
\end{document}